\documentclass[11pt]{article}
\usepackage{amsfonts}
\usepackage{graphicx}

\newcommand{\e}{\varepsilon}
\newcommand{\tQ}{\tilde Q}
\newcommand{\tI}{\tilde I}
\newcommand{\R}{\mathbb R}
\newcommand{\proof}{{\bf Proof:}\\}
\newcommand{\eproof}{{$\square$}\\}

\newtheorem{lemma}{Lemma}[section]
\newtheorem{proposition}{Proposition}[section]
\newtheorem{theorem}{Theorem}[section]
\newtheorem{remark}{Remark}

\title{Fast conservative and entropic numerical methods for the Boson Boltzmann
equation\thanks{This work was supported by the WITTGENSTEIN AWARD
2000 of Peter Markowich, financed by the Austrian Research Fund
FWF and by the European network HYKE, funded by the EC as contract
HPRN-CT-2002-00282.}}
\author{Peter A.~Markowich\thanks{University of Vienna, Institute of Mathematics,
Boltzmanngasse 9 A - 1090 Vienna, Austria ({\tt
Peter.Markowich@univie.ac.at})} \and Lorenzo
Pareschi\thanks{University of Ferrara, Department of Mathematics,
Via Machiavelli 35 - 44100 Ferrara, Italy ({\tt
lorenzo.pareschi@unife.it})}}
\begin{document}

\date{November 4, 2004}
\maketitle

\begin{abstract}
In this paper we derive accurate numerical methods for the quantum
Boltzmann equation for a gas of interacting bosons. The schemes
preserve the main physical features of the continuous problem,
namely conservation of mass and energy, the entropy inequality and
generalized Bose-Einstein distributions as steady states. These
properties are essential in order to develop schemes that are able
to capture the energy concentration behavior of bosons. In
addition we develop fast algorithms for the numerical evaluation
of the resulting quadrature formulas which allow the final schemes
to be computed only in $O(N^2\log_2 N)$ operations instead of
$O(N^3)$.
\end{abstract}
\bigskip
\begin{bf}
Key words:
\end{bf}
Boson Boltzmann equation, condensation, quadrature formulas, fast
algorithms.
\\
\\
\begin{bf}
AMS Subject classification:
\end{bf}
82C10, 76P05, 65D32, 65T50


\section{Introduction}
We consider a gas of interacting bosons, which are trapped by a
confining potential $V=V(x)$ with $\min V(x)=0$. We denote the
total energy of a boson with momentum $p$ and position $x$ (after
an appropriate non-dimensionalization) by
\begin{equation}
\e(x,p)=\frac{|p|^2}{2}+V(x).
\end{equation}
Let $F=F(p,x,t)\geq 0$ be the phase-space density of bosons.
Assuming a boson distribution which only depends on the total
energy $\e$ we write
\begin{equation}
F(x,p,t)=f\left(\frac{|p|^2}{2}+V(x),t\right),
\label{eq1:psd}
\end{equation}
where $f=f(\e,t)\geq 0$ is the boson density in energy space.

\subsection{The Boson Boltzmann equation}
Following
\cite{ST1},\cite{ST2},\cite{Jaksh1},\cite{GZ1},\cite{GZ2} we write
a Boltzmann-type equation (referred to as boson Boltzmann equation
in the sequel) in energy space
\begin{equation}
\rho(\e)\frac{\partial f}{\partial t}= Q(f)(\e), \quad t >0,
\label{eq1:qbe}
\end{equation}
with the collision integral
\begin{eqnarray}
\nonumber
Q(f)(\e)&=&\int_{\R_+^3}\delta(\e+\e_*-\e'-\e_*')S(\e,\e_*,\e',\e_*')[f'
f_*' (1+f)(1+f_*)\\
\label{eqn:qbe}
\\[-.2cm]\nonumber
&-& f f_* (1+f')(1+f_*')]\,d\e_*d\e'd\e_*',
\end{eqnarray}
where $S\geq 0$ is a given function.

We denoted the density of states by
\begin{equation}
\rho(\e)=\int_{\R^6}\delta\left(\e-\left(\frac{|p|^2}{2}+V(x)\right)\right)dp\,dx,
\label{eq1:ds}
\end{equation}
and
\begin{equation}
f'=f(\e',t),\quad f_*'=f(\e_*',t),\quad f=f(\e,t),\quad
f_*=f(\e_*,t).
\end{equation}
As usual $\e$ and $\e_*$ are the pre-collisional energies of two
interacting bosons and $\e'$ and $\e'_*$ are the post-collisional
ones.

The positive measure
\begin{equation}
\delta(\e+\e_*-\e'-\e_*')S(\e,\e_*,\e',\e_*')
\label{eq1:etr}
\end{equation}
denotes the energy transition rate, i.e. $Sd\e'\,d\e'_*$ is the
transition probability per unit volume and per unit time that two
bosons with incoming energies $\e$, $\e_*$ are scattered with
outgoing energies $\e'$, $\e'_*$.

A simple computation shows that the phase-space density
$F=F(x,p,t)$ satisfies the momentum-position space Boltzmann
equation
\begin{equation}
\frac{\partial F}{\partial t}+p\cdot \nabla_x F-\nabla_x V(x)\cdot
\nabla_p F={\tilde Q}(F),
\label{eq1:mpqbe}
\end{equation}
with the scattering integral
\begin{equation}
{\tilde
Q}(F)(x,p)=\frac{Q(F)\left({|p|^2}/{2}+V(x)\right)}{\rho\left({|p|^2}/{2}+V(x)\right)}.
\label{eq:SI}
\end{equation}
Note that (\ref{eq:SI}) does not correspond to the physical
Boltzmann operator for bosons except in the homogeneous case
$V(x)=0$ and $F$ independent of $x$, where we set
\begin{equation}
\rho(\e)=\int_{\R^3}\delta\left(\e-\frac{|p|^2}{2}\right)dp
\label{eq1:dsh}
\end{equation}
and compute
\begin{equation}
\rho(\e)=4\pi \sqrt{2\e}. \label{eq1:dshc}
\end{equation}
Then equation (\ref{eq1:mpqbe}) is formally identical to the Boson
Boltzmann equation considered in \cite{EMV1},\cite{EMV2}
\begin{eqnarray}
\nonumber
\frac{\partial F}{\partial
t}&=&\int_{\R^9}\delta(p+p_*-p'-p_*')\delta(\e+\e_*-\e'-\e_*')
W(p,p_*,p',p_*')\\
\label{eq1:mpqbe2}
\\[-.25cm]
\nonumber &&[F'F_*' (1+F)(1+F_*)- F F_*
(1+F')(1+F_*')]\,dp_*dp'dp_*',
\end{eqnarray}
with $\e(p)=|p|^2/2$ and $W$, $S$ are related by
\begin{eqnarray}
\nonumber &&\int_{S^2\times S^2\times S^2 }\delta(p+p_*-p'-p_*')
W(p,p_*,p',p_*')\,d\sigma_*d\sigma'd\sigma_*'\\
\nonumber
&=&\frac{S\left({|p|^2}/{2},{|p_*|^2}/{2},{|p'|^2}/{2},{|p_*'|^2}/{2}\right)}
{\rho({|p|^2}/{2})|p_*||p'||p'_*|}.
\end{eqnarray}
Here we denoted $p_*=|p_*|\sigma_*$, $p'=|p'|\sigma'$,
$p=|p|\sigma$, and $p_*'=|p_*'|\sigma_*'$. In particular for
$W\equiv 1$ we have
\begin{equation}
S(\e,\e_*,\e',\e_*')={\rm const}\,\rho(\e_{\min}),
\end{equation}
where (see \cite{EMV1})
\begin{equation}
\e_{\min}=\min(\e,\e_*,\e',\e_*').
\end{equation}
Even in the non-homogeneous case $V(x)\neq 0$ the equation
(\ref{eq1:qbe}) is formally identical to the isotropic version of
the homogeneous bosonic Boltzmann equation (\ref{eq1:mpqbe2})
(after the introduction of $|p|^2/2$ as new independent variable).
However, the density of states is computed by formula
(\ref{eq1:ds}) in the non homogeneous case instead of
(\ref{eq1:dsh}) in the space homogeneous case.

In the physical literature the equation (\ref{eq1:qbe}), usually
referred to as {\em ergodic approximation} of the Boltzmann
equation, is derived in the nonhomogeneous case as approximation
of the phase-space Boltzmann equation by a projection technique
\cite{Jaksh1},\cite{LRW}.

For a mathematical analysis of the bosonic Boltzmann equation in
the space homogeneous isotropic case we refer to
\cite{Lu1},\cite{Lu2},\cite{EMV1},\cite{EMV2}. We remark that
already the issue of giving mathematical sense to the collision
operator $Q(f)$ is highly nontrivial (particularly for scattering
rates without cutoff or if positive measures $f$ are allowed, as
required by a careful analysis of the equilibrium states).

\subsection{Physical properties}

A simple calculation gives the weak form of the collision
operator. Let $\phi=\phi(\e)$ be a test function. Then, at least
formally
\begin{eqnarray}
\nonumber \int_0^\infty Q(f)\phi d\e&=&\frac12
\int_{\R_+^4}\delta(\e+\e_*-\e'-\e_*')S(\e,\e_*,\e',\e_*')[f'
f_*' (1+f)(1+f_*)\\
\label{eq1:qbewc}
\\[-.2cm]\nonumber
&-& f f_* (1+f')(1+f_*')][\phi+\phi_*-\phi'-\phi_*']d\e
d\e_*d\e'd\e_*'.
\end{eqnarray}
Here we used the micro-reversibility property, i.e. the fact that
each collision is reversible and that each pair of interacting
bosons represents a closed physical system. Mathematically this
amounts to the requirement \cite{EMV1}
\begin{equation}
S(\e,\e_*,\e',\e_*')=S(\e_*,\e,\e',\e_*')=S(\e',\e_*',\e,\e_*).
\label{eq1:sp}
\end{equation}
The symmetry properties (\ref{eq1:sp}) immediately imply the
analogous properties for the energy transition rate
(\ref{eq1:etr}) and the weak form (\ref{eq1:qbewc}) follows from
the variable substitution in the integral using these symmetries.

As a consequence we have the following collision invariants
\begin{enumerate}
\item
\begin{equation}
\phi(\e)\equiv 1\quad \Rightarrow \quad \int_{0}^\infty
Q(f)(\e)\,d\e=0, \label{eq1:mc}
\end{equation}
\item
\begin{equation}
\phi(\e)\equiv \e\quad \Rightarrow \quad \int_{0}^\infty
Q(f)(\e)\e\,d\e=0. \label{eq1:ec}
\end{equation}
\end{enumerate}

Consider now the IVP (\ref{eq1:qbe}) supplemented by the initial
condition
\begin{equation}
f(\e,t=0)=f_0(\e)\geq 0,\quad \e>0. \label{eq1:qbei}
\end{equation}
Then (\ref{eq1:mc}) implies mass conservation
\begin{equation}
\int_{0}^\infty \rho(\e)f(\e,t)\,d\e=\int_{0}^\infty
\rho(\e)f_0(\e)\,d\e,\quad \forall\,t>0, \label{eq:N}
\end{equation}
and (\ref{eq1:ec}) energy conservation
\begin{equation}
\int_{0}^\infty \rho(\e)f(\e,t)\e\,d\e=\int_{0}^\infty
\rho(\e)f_0(\e)\e\,d\e,\quad \forall\,t>0. \label{eq:E}
\end{equation}

The H-theorem for (\ref{eq1:qbe}) is derived by setting
$\phi(\e)=\ln(1+f(\e))-\ln f(\e)$ in (\ref{eq1:qbewc}).

We calculate
\begin{eqnarray}
\nonumber &&\int_0^\infty Q(f)(\e)(\ln(1+f(\e))-\ln
f(\e))d\e\\
\label{eq1:qbewch}
\\[-.2cm]\nonumber
&=&
\frac12\int_{\R_+^4}\delta(\e+\e_*-\e'-\e_*')S(\e,\e_*,\e',\e_*')
e(f)d\e d\e_*d\e'd\e_*':=D[f],
\end{eqnarray}
where
\begin{equation}
e(f)=z(ff_*(1+f')(1+f'_*),f'f'_*(1+f)(1+f_*))
\end{equation}
and
\begin{equation}
z(x,y)=(x-y)(\ln x - \ln y).
\end{equation}
Since the integrand of the entropy dissipation $D[f]$ is
non-negative, we deduce the following H-theorem, obtained by
multiplying (\ref{eq1:qbe}) by $\phi(\e)=\ln(1+f(\e))-\ln f(\e)$
\begin{equation}
\frac{d}{dt}S[f]=D[f],
\end{equation}
which implies that the entropy
\begin{equation}
S[f]:=\int_0^\infty \rho(\e)((1+f)\ln(1+f)-f\ln f)d\e,
\end{equation}
is increasing along trajectories of (\ref{eq1:qbe}). We remark
that trivially the third physical conservation law, namely
momentum conservation, also holds. Clearly the phase-space density
$F$ of (\ref{eq1:psd}) satisfies
\begin{equation}
\int_{\R^3} p F(x,p,t)dx\equiv 0, \quad \forall t \geq 0.
\end{equation}
We now turn to the issue of steady states of (\ref{eq1:qbe}). The
problem of equilibrium distributions for bosons has a very long
history, going back to Bose and Einstein in the twenties of the
last century (see \cite{Bo},\cite{E1},\cite{E2}), who noticed that
the class of 'regular' Bose-Einstein distributions
\begin{equation}
f_\infty(\e)=\frac{1}{e^{\alpha\e+\beta}-1},\quad \alpha >0,
\beta>0 \label{eq:RBE}
\end{equation}
is not sufficient to assume all arbitrarily large values of
equilibrium mass
\begin{equation}
M_\infty=\int_0^\infty \rho(\e) f_\infty(\e)d\e, \label{eq1:em}
\end{equation}
and arbitrarily small values of equilibrium energy
\begin{equation}
E_\infty=\int_0^\infty \rho(\e)\e f_\infty(\e) d\e, \label{eq1:ee}
\end{equation}
such that Dirac distributions centered in zero energy have to be
included in the set of equilibrium states. In \cite{EMV1} it was
shown that for every pair $(M_\infty,E_\infty)\in \R^2_+$ there
exist $\alpha \geq 0$, $\beta \in \R$ such that the generalized
Bose-Einstein distribution defined by
\begin{equation}
{\rho(\e)}f_\infty(\e)=\frac{{\rho(\e)}}{e^{\alpha\e+\beta_+}-1}+{|\beta_-|}\delta(\e),
\label{eq:GBE}
\end{equation}
is an equilibrium state of (\ref{eq1:qbe}) (in the sense of
maximizing the entropy, see \cite{EMV1} for analytical details)
satisfying (\ref{eq1:em})-(\ref{eq1:ee}). Here we denoted
$\beta_+=\max(\beta,0)$ and $\beta_-=-\max(-\beta,0)$. The value
$M_{\infty,cond}=|\beta_-|$ represents the mass of particles which
are condensed in equilibrium, i.e. in their quantum mechanical
ground state with $\e=0$.

Off course, it is analytically nontrivial to define the
nonlinearities in the entropy (dissipation) and in the collision
operator, in particular for measures which are singular with
respect to the Lebesgue measure, as required for the equilibrium
states. For details we refer to the references \cite{EMV1} and
\cite{Lu1}, here we only mention that an approximation argument
shows that the singular part of a measure $f$  does not contribute
to the entropy $S[f]$. For appropriate scattering rates (with
unphysical cut-off) in the homogeneous case an
existence/uniqueness theory for integrable and for measure
solutions can be set up. So far, it is not clear how the cut-off
assumption can be removed.

In the following sections we shall use
\begin{equation}
S(\e,\e_*,\e',\e'_*)=\rho(\e_{\min}). \label{eq1:ass}
\end{equation}

Notice that the condensation is fully localized in phase space,
i.e. it may only occur at $p=0$ (vanishing momentum) and at those
points in position space, where the potential assumes its minimum
value $0$. The reason for this is the form (\ref{eq1:psd}) of the
phase space distribution and a semiclassical limit process which
leads to the Boson Boltzmann equation (\ref{eq1:qbe}).

The purpose of this paper is to derive an accurate discretization
of the IVP (\ref{eq1:qbe}), (\ref{eq1:qbei}), which maintains the
basic analytical and physical features of the continuous problem,
namely
\begin{itemize}
\item Mass and energy conservation
\item Entropy growth
\item Generalized Bose-Einstein equilibrium distribution
\end{itemize}
To this aim we shall derive first and second order accurate
quadrature formulas for $Q(f)$. These schemes due to their
'direct' derivation from the continuous operator possess all the
desired physical properties at a discrete level. In addition we
show that with the choice (\ref{eq1:ass}) the computations can be
performed with a fast algorithms reducing the $O(N^3)$ cubic cost
to $O(N^2\log_2 N)$. For the sake of completeness we mention the
recent works
\cite{bcdl},\cite{Lemou},\cite{Pa},\cite{PRT2},\cite{PTV} in which
fast methods for Boltzmann equations were derived using different
techniques like multipole methods, multigrid methods and spectral
methods.

The rest of the paper is organized as follows. In the next Section
we discuss the details of our numerical schemes, together with the
issues of consistency and computational complexity. In Section 3
several numerical tests are performed. The results confirm the
expected accuracy of the schemes and in particular show the
ability of the methods to capture the concentration behavior of
bosons. Finally we concluded the paper with some remarks in
Section 4.

\section{Fast, conservative and entropic methods}
We consider the IVP for the quantum boson Boltzmann equation
\begin{eqnarray}
\label{eqn:qbo}
\rho(\e)\frac{\partial f}{\partial t}&=& Q(f)(\e), \quad t >0,\\
\label{eqn:qboi}
f(\e,t=0)&=&f_0(\e)\geq 0.
\end{eqnarray}
Here the independent variable $\e>0$ represents the kinetic
energy, $\rho=\rho(\e)\geq 0$ is the (given) density of states and
the boson collision operator now reads
\begin{eqnarray}
\nonumber
Q(f)(\e)&=&\int_{\R_+^3}\delta(\e+\e_*-\e'-\e_*')\rho(\e_{\min})[f'
f_*' (1+f)(1+f_*)\\
\label{eqn:qbes}
\\[-.2cm]\nonumber
&-& f f_* (1+f')(1+f_*')]\,d\e_*d\e'd\e_*'.
\end{eqnarray}
Obviously the equation (\ref{eqn:qbo}) maintains a minimum
principle such that solution of (\ref{eqn:qbo}), (\ref{eqn:qboi})
satisfy $f(\e,t)\geq0$ for $\e\geq0, t>0$ if $f_0(\e)\geq0$ for
$\e>0$.

\subsection{Reduction on a bounded domain}

Our starting point in the development of a numerical scheme for
(\ref{eqn:qbes}) is the definition of a bounded domain
approximation of the collision operator $Q$.

Let $f$ be defined for $\e \in [0,R]$ and denote
\begin{eqnarray}
\nonumber
Q_R(f)(\e)&=&\int_{[0,R]^3}\delta(\e+\e_*-\e'-\e_*')\rho(\e_{\min})[f'
f_*' (1+f)(1+f_*)\\
\label{eqn:qbec}
\\[-.2cm]\nonumber
&-& f f_* (1+f')(1+f_*')]\psi(\e \leq R)\,d\e_*d\e'd\e_*'
\end{eqnarray}
where $\psi(I)$ is the indicator function of the set $I$. Then, at
least formally
\begin{eqnarray}
\nonumber \int_0^\infty Q_R(f)\phi d\e&=&\frac12
\int_{[0,R]^4}\delta(\e+\e_*-\e'-\e_*')\rho(\e_{\min})[f'
f_*' (1+f)(1+f_*)\\
\label{eqn:qbewc}
\\[-.2cm]\nonumber
&-& f f_* (1+f')(1+f_*')][\phi+\phi_*-\phi'-\phi_*']d\e
d\e_*d\e'd\e_*'
\end{eqnarray}
for any test function $\phi=\phi(\e)$.
The proof follows the lines of the corresponding weak form of $Q$
discussed in Section 1.

Consider now the approximate IVP
\begin{eqnarray}
\label{eqn:qbor}
\rho(\e)\frac{\partial f_R}{\partial t}&=& Q_R(f_R)(\e), \quad t >0,\\
\label{eqn:qboir} f_R(\e,t=0)&=&f_{0,R}(\e)\geq 0.
\end{eqnarray}

Then the weak formulation (\ref{eqn:qbewc}) of $Q_R$ implies mass
and energy conservation
\begin{equation}
\int_{0}^R \rho(\e)f_R(\e,t)\,d\e=\int_{0}^R
\rho(\e)f_{0,R}(\e)\,d\e,\quad \forall\,t>0,
\end{equation}
\begin{equation}
\int_{0}^R \rho(\e)f_R(\e,t)\e\,d\e=\int_{0}^R
\rho(\e)f_{0,R}(\e)\e\,d\e,\quad \forall\,t>0.
\end{equation}
Also the entropy inequality
\begin{equation}
\frac{d}{dt}S_R[f_R]\geq 0,
\end{equation}
holds with the entropy
\begin{equation}
S_R[f_R]:=\int_0^R \rho(\e)((1+f_R)\ln(1+f_R)-f_R\ln f_R)d\e.
\end{equation}
These properties are in full analogy with the corresponding ones
of the IVP (\ref{eq1:qbe}), (\ref{eq1:qbei}).

\subsection{Discretization and main properties}

Let us now introduce the set of discrete energy grid points
$\e_1\leq \e_2\leq \ldots \leq \e_N$ in $[0,R]$. A general
quadrature formula for (\ref{eqn:qbec}) is given by

\begin{eqnarray}
\nonumber Q_R(f)(\e_i)\approx {\tQ}_R(f)(\e_i)&=&\sum_{j,k,l=1}^N
w_{ij}^{kl}\,\delta_{ij}^{kl}\, \rho(\e_{\min}) [f_k
f_l(1+f_i)(1+f_j)\\
\label{eqn:qbed}
\\[-.2cm]\nonumber
&-&f_i f_j (1+f_k)(1+f_l)]\psi(\e_i \leq R),
\end{eqnarray}
where now $f_i=f(\e_i)$ and
$\e_{\min}=\min\{\e_i,\e_j,\e_k,\e_l\}$. The quantities
$w_{ij}^{kl}$ are the weights of the quadrature formula and
$\delta_{ij}^{kl}$ a suitable discretization of the
$\delta$-function on the grid.

In order to maintain the conservation properties on the discrete
level it is of paramount importance that the discretized
$\delta$-function will reduce the points in the sum to a discrete
index set which satisfies the relation $i+j=k+l$. Thus it is
natural to restrict to equally spaced grid points which satisfy
exactly the aforementioned relation on the computational grid.

We will further simplify the quadrature formula by considering
product quadrature rules with equal weights for which
$w_{ij}^{kl}=w_j w_k w_l=w^3$ with $w=R/N$ and
\[
\int_0^R f(\e)\,d\e \approx w \sum_{i=1}^N f(\e_i).
\]

We now consider the set of ODEs which originates from the energy
discretization of the IVP (\ref{eqn:qbor}), (\ref{eqn:qboir})
\begin{eqnarray}
\label{eqn:qbord}
\rho(\e_i)\frac{d f_i}{d t}&=& \tQ_R(f)(\e_i), \quad t >0,\\
\label{eqn:qboird} f_i(t=0)&=&f_{0,R}(\e_i)\geq 0.
\end{eqnarray}
and prove
\begin{proposition}
If we define
\begin{equation}
\delta_{ij}^{kl}=\left\{
\begin{array}{cc}
  1/w & i+j=k+l \\
  0 & {\rm otherwise}
\end{array}
\right.
\end{equation}
the solutions of the IVP (\ref{eqn:qbord}), (\ref{eqn:qboird})
satisfy the following discrete conservation properties and entropy
principle
\begin{equation}
w \sum_{i=1}^N \rho(\e_i) \frac{d f_i}{d t} \phi(\e_i)=0,\qquad
\phi(\e)=1,\quad \phi(\e)=\e, \label{eq:consd}
\end{equation}
\begin{equation}
w \sum_{i=1}^N \rho(\e_i) \frac{d h(f_i)}{d t} \geq 0,\qquad
h(f_i)=(1+f_i)\log (1+f_i)-f_i\log f_i. \label{eq:entd}
\end{equation}
\end{proposition}
\proof

Due to the definition of $\delta_{ij}^{kl}$ we have the quadrature
formula
\begin{eqnarray}
\nonumber {\tQ}_R(f)(\e_i)&=&w^2 \sum_{{j,l=1}\atop{1\leq k=i+j-l
\leq N}}^N\, \rho(\e_{\min}) [f_k
f_l(1+f_i)(1+f_j)\\
\label{eqn:qbed1}
\\[-.2cm]\nonumber
&-&f_i f_j (1+f_k)(1+f_l)].
\end{eqnarray}

In particular, for any test function $\phi$, formula
(\ref{eqn:qbed1}) admits the following discrete analogous of the
corresponding weak identity for the collision operator
\begin{eqnarray}
\nonumber w \sum_{i=1}^N {\tQ}_R(f)(\e_i)\phi(\e_i) &=&\frac12 w^3
\sum_{{i,j,k,l=1}\atop{i+j=k+l}}^N\, \rho(\e_{\min}) [f_k
f_l(1+f_i)(1+f_j)\\
\label{eqn:qbedw}
\\[-.2cm]\nonumber
&-&f_i f_j (1+f_k)(1+f_l)][\phi_i+\phi_j-\phi_k-\phi_l],
\end{eqnarray}
where $\phi_i=\phi(\e_i)$. The equations (\ref{eq:consd}) are
obtained taking $\phi(\e)=1$, and $\phi(\e)=\e$. The discrete
entropy inequality can be derived choosing
$\phi(\e)=h'(f(\e))=\ln(1+f(\e))-\ln f(\e)$. In fact, as in the
continuous case, we find
\begin{eqnarray}
\nonumber &&w \sum_{i=1}^N \rho(\e_i) \frac{d h(f_i)}{d t}=
\frac12 w^3 \sum_{{i,j,k,l=1}\atop{i+j=k+l}}^N\, \rho(\e_{\min})
[f_k
f_l(1+f_i)(1+f_j)\\
\label{eq:weekd}
\\[-.2cm]\nonumber
&-&f_i f_j (1+f_k)(1+f_l)][h'(f_i)+h'(f_j)-h'(f_k)-h'(f_l)] \geq
0,
\end{eqnarray}
since
\[
h'(f_i)+h'(f_j)-h'(f_k)-h'(f_l)=\log((1+f_i)(1+f_j)f_k f_l)-
\log((1+f_k)(1+f_l)f_i f_j),
\]
and the function $z(x,y)=(x-y)(\log x - \log y) \geq 0$ for $x,y
\in \R^+$.

 \eproof

\begin{remark}
It is easy to check by direct verification using (\ref{eq:weekd})
that these schemes admits 'regular' discrete Bose-Einstein
equilibrium states of the form
\begin{equation}
f_\infty(\e_i)=\frac{1}{e^{\alpha\e_i+\beta}-1},\quad \alpha >0,
\beta \in \R. \label{eq:BEn}
\end{equation}
More delicate is the question of 'generalized' discrete
Bose-Einstein equilibrium which will be discussed later on.
\end{remark}

\begin{remark}
Clearly one may use other product quadrature rules with different
weights. However then the definition of a consistent discrete
$\delta$-function which satisfies the aforementioned conservation
laws and entropy principle becomes very difficult. On the other
hand it is shown in the next section that the choice of quadrature
(\ref{eqn:qbed1}) includes numerical methods up to second order
accuracy.
\end{remark}

\subsection{First and second order methods}

\begin{figure}[htb]
\centerline{
\includegraphics[scale=0.6]{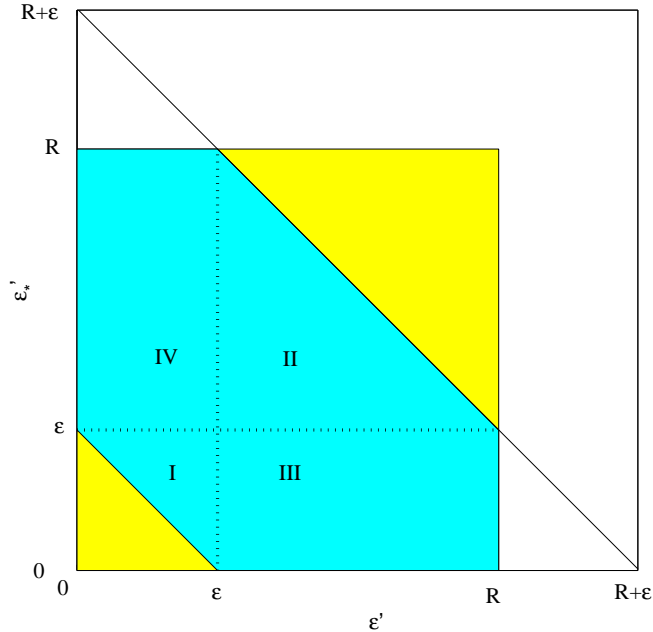}
} \caption{ \small The computational domain (dark gray region) in
the $(\epsilon', \epsilon_*')$ plane for a fixed $\e$}
\label{fig:domain}
\end{figure}

Let us rewrite for $\e \in [0,R]$ the collision integral
(\ref{eqn:qbec}) as

\begin{equation}
Q_R(f)(\e)=\int_0^R \int_{S(\e,\e')}^{D(\e,\e')} \rho(\e_{\min})
F(\e,\e',\e'_*)\,d\e_*'d\e', \label{eqn:qbec2}
\end{equation}
where $F(\e,\e',\e'_*)=[f'f_*' (1+f)(1+f_*)- f f_*
(1+f')(1+f_*')]$, with $\e_*=\e'+\e_*'-\e$, and
$S(\e,\e')=\max\{\e-\e',0\}$, $D(\e,\e')=\min\{\e-\e'+R,R\}$. The
integration domain for a fixed value of $\e$ in the $(\e',\e_*')$
plane is shown in figure \ref{fig:domain}.

We need the following

\begin{lemma}
We have
\begin{equation}
\rho(\e_{\min})=\left\{
\begin{array}{cc}
  \rho(\e_*) & \quad (\e',\e_*')\in {\rm I}\\
  \rho(\e) &  \quad (\e',\e_*')\in {\rm II}\\
  \rho(\e_*') &  \quad (\e',\e_*')\in {\rm III}\\
  \rho(\e') &  \quad (\e',\e_*')\in {\rm IV}
\end{array}
\right.
\end{equation}
where the regions I, II, III, IV represent a partition of the
computational domain and are shown in figure \ref{fig:domain}.
\end{lemma}
\proof

Region I is characterized by $0 \leq \e_*'\leq \e$ and $0 \leq \e'
\leq \e$ with $\e'+\e'_* \geq \e$. Thus $\e_*=\e_*'+\e' -\e \leq
\e$, $\e_*=\e_*'+\e' -\e \leq \e'$, $\e_*=\e_*'+\e' -\e \leq
\e_*'$ and hence $\e_{\min}=\e_*$.

Region II is characterized by $\e'\geq \e$ and $\e'_* \geq \e$
with $\e'+\e'_* \leq R+\e$. Thus $\e_*=\e'+\e'_* -\e \geq \e$ and
hence $\e_{\min}=\e$.

Region III is characterized by $R\geq \e'\geq \e$ and $0\leq \e'_*
\leq \e$. Thus $\e_*=\e'+\e'_* -\e \geq \e'_*$ and hence
$\e_{\min}=\e'_*$.

Region IV is characterized by $R\geq \e_*'\geq \e$ and $0\leq \e'
\leq \e$. Thus $\e_*=\e'+\e'_* -\e \geq \e'$ and hence
$\e_{\min}=\e'$.

\eproof

Using the previous lemma the integral (\ref{eqn:qbec2}) over the
four regions can be decomposed as
\begin{equation}
Q_R(f)(\e)=I_1(\e)+I_2(\e)+I_3(\e)+I_4(\e),
\end{equation}
with
\begin{eqnarray}
I_1(\e)&=&\int_0^\e \int_{\e-\e'}^{\e} \rho(\e'+\e_*'-\e)
F(\e,\e',\e'_*)\,d\e_*'d\e',\\
I_2(\e)&=&\int_\e^R \int_{\e}^{R+\e-\e'} \rho(\e)
F(\e,\e',\e'_*)\,d\e_*'d\e',\\
I_3(\e)&=&\int_\e^R \int_{0}^{\e} \rho(\e_*')
F(\e,\e',\e'_*)\,d\e_*'d\e',\\
I_4(\e)&=&\int_0^\e \int_{\e}^{R} \rho(\e')
F(\e,\e',\e'_*)\,d\e_*'d\e'.
\end{eqnarray}

A similar decomposition holds for the quadrature formula
(\ref{eqn:qbed1})
\begin{equation}
{\tQ}_R(f)(\e_i)={\tI}_1(\e_i)+{\tI}_2(\e_i)+{\tI}_3(\e_i)+{\tI}_4(\e_i),
\label{eqn:qbed1d}
\end{equation}
with
\begin{eqnarray}
{\tI}_1(\e_i)&=&w^2 \sum_{k=1}^i\sum_{l=i-k+1}^{i}\,
\rho(\e_k+\e_l-\e_i)
F(\e_i,\e_k,\e_l),\\
{\tI}_2(\e_i)&=&w^2 \sum_{k=i+1}^N\sum_{l=i}^{N+i-k}\, \rho(\e_i)
F(\e_i,\e_k,\e_l),\\
{\tI}_3(\e_i)&=&w^2 \sum_{k=i+1}^N\sum_{l=1}^{i}\, \rho(\e_l)
F(\e_i,\e_k,\e_l),\\
{\tI}_4(\e_i)&=&w^2 \sum_{k=1}^i\sum_{l=i+1}^{N}\, \rho(\e_k)
F(\e_i,\e_k,\e_l).
\end{eqnarray}

 From the point of view of accuracy we can state

\begin{theorem}[Consistency]
Let the function $f$ and $\rho$ be $C^m([0,R])$, $m=1$ or $m=2$,
then the quadrature formula (\ref{eqn:qbed1}) satisfies
\begin{equation}
|Q_R(f)(\e_i)-{\tQ}_R(f)(\e_i)| \leq R^2
C_m(\Delta\e)^{m}M_m,\qquad \Delta\e=R/N,
\end{equation}
where $M_m$ is a constant that depends on $f$ and $\rho$ and their
derivatives up to the order $m$ and if $\e_i=(i-1)\Delta\e$,
$i=1,\ldots, N$ {\em (rectangular rule)} then $m=1$ and $C_m=1/2$,
whereas if $\e_i=(i-1/2)\Delta\e$, $i=1,\ldots, N$ {\em (midpoint
rule)} $m=2$ and $C_m=1/24$.
\end{theorem}
\proof First let us recall the following basic estimate for a
composite product quadrature rule with equal weights (see
\cite{Rab} for example)
\begin{eqnarray}
\nonumber
\left| \int_a^b \int_c^d g(x,y)dxdy-\Delta x \Delta y
\sum_{i=1}^{N_x}\sum_{j=1}^{N_y} g(x_i,y_j)\right| \leq \\
\label{eqn:ees}
\\[-.25cm]
\nonumber {(b-a)(d-c)}C_m\left[(\Delta x)^m M_{x,m}+ (\Delta y)^m
M_{y,m}\right],
\end{eqnarray}
where $\Delta x=(b-a)/N_x$, $\Delta y=(d-c)/N_y$, $M_{x,m}$ and
$M_{y,m}$ are two constants such that
\[
|\frac{\partial^m g}{\partial x^m}| \leq M_{x,m}, \quad
|\frac{\partial^m g}{\partial y^m}| \leq M_{y,m},
\]
on $[a,b]\times[c,d]$ and if $x_i=(i-1)\Delta x$, $y_i=(i-1)\Delta
y$ then $m=1$, $C_m=1/2$, whereas if $x_i=(i-1/2)\Delta x$,
$y_i=(i-1/2)\Delta y$ then $m=2$, $C_m=1/24$.

Now, since the integrands which appear in $I_i$ satisfy the
required regularity conditions and approximations given by
${\tI}_i$ are the corresponding generalized composite product
quadrature rules, each error $|I_i-{\tI}_i|$ can be estimated
similarly to (\ref{eqn:ees}). More precisely we have
\begin{eqnarray*}
|I_1(\e_i)-{\tI}_1(\e_i)|&\leq& {(\e_i)^2}C_m(\Delta
\e)^m[M^1_{\e',m}(\e_i)+
M^1_{\e_*',m}(\e_i)],\\
|I_2(\e_i)-{\tI}_2(\e_i)|&\leq& {(R-\e_i)^2}C_m(\Delta
\e)^m[M^2_{\e',m}(\e_i)+
M^2_{\e_*',m}(\e_i)],\\
|I_3(\e_i)-{\tI}_3(\e_i)|&\leq& {\e_i(R-\e_i)}C_m(\Delta
\e)^m[M^3_{\e',m}(\e_i)+
M^3_{\e_*',m}(\e_i)],\\
|I_4(\e_i)-{\tI}_4(\e_i)|&\leq& {\e_i(R-\e_i)}C_m(\Delta
\e)^m[M^4_{\e',m}(\e_i)+ M^4_{\e_*',m}(\e_i)],
\end{eqnarray*}
where the constants $M^i_{\e',m}(\e)$ and $M^i_{\e_*',m}(\e)$ are
suitable bounds of the partial derivatives of order $m$ of the
integrand functions.

Summing up the errors we get
\begin{equation}
|Q_R(f,f)(\e_i)-{\tQ}_R(\e_i)|\leq {R^2}C_m(\Delta \e)^m M_m,
\end{equation}
where
$M_m(\e)=\max_{i,k}\{M^i_{\e',m}(\e_k)+M^i_{\e_*',m}(\e_k)\}$.

\eproof

\subsection{Fast algorithms}
\label{fast}

Finally we will analyze the problem of the computational cost of
the quadrature formula (\ref{eqn:qbed1}). A straightforward
analysis shows that the evaluation of the double sum in
(\ref{eqn:qbed1}) at the point $\e_i$ requires
$(2(i-1)(N-i+1)+N^2)/2$ operations. The overall cost for all $N$
points is then approximatively $2N^3/3$. However using transform
techniques and the decomposition (\ref{eqn:qbed1d}) this $O(N^3)$
cost can be reduced to $O(N^2\log_2 N)$.

In order to do this let us set $h=k+l=i+j$ in (\ref{eqn:qbed1})
and rewrite
\begin{eqnarray}
\nonumber {\tQ}_R(\e_i)&=&w^2 \sum_{h=2}^{2N} \sum_{k=1}^N\,
\rho(\e_{\min}) [f_{k}
f_{h-k}(1+f_i)(1+f_{h-i})\\
\label{eqn:qbed1f}
\\[-.2cm]\nonumber
&-&f_i f_{h-i}
(1+f_{k})(1+f_{h-k})]\Psi_{h-i}^{[1,N]}\Psi_{h-k}^{[1,N]},
\end{eqnarray}
where we have set
\begin{equation}
\Psi_{i}^{[s,d]}=\left\{
\begin{array}{cc}
  1 & s \leq i \leq d \\
  0 & {\rm otherwise}
\end{array}
\right.
\end{equation}
In (\ref{eqn:qbed1f}) we assume that the function $f_i$ is
extended to $i=1,\ldots,2N$ by padding zeros for $i > N$.

The sum (\ref{eqn:qbed1f}) can be split into sum over the four
regions which characterize $\rho(\e_{\min})$. We shall give the
details of the fast algorithm only for region I, the other regions
can be treated similarly. We have
\begin{eqnarray}
\nonumber {\tI}_1(\e_i)&=&w^2 \sum_{h=2}^{2N}\sum_{k=1}^{i}\,
\rho(\e_{h-i})[f_{k}
f_{h-k}(1+f_i)(1+f_{h-i})\\
\label{eqn:I1f}
\\[-.2cm]\nonumber
&-&f_i f_{h-i}
(1+f_{k})(1+f_{h-k})]\Psi_{h-i}^{[1,i]}\Psi_{h-k}^{[1,i]},
\end{eqnarray}
or equivalently
\begin{eqnarray*}
\nonumber {\tI}_1(\e_i)&=&w^2
\sum_{h=2}^{2N}\rho(\e_{h-i})(1+f_i)(1+f_{h-i})\Psi_{h-i}^{[1,i]}
\sum_{k=1}^{i}\,f_{k} f_{h-k}\Psi_{h-k}^{[1,i]}\\
&-& w^2\sum_{h=2}^{2N}\rho(\e_{h-i}) f_i f_{h-i}
\Psi_{h-i}^{[1,i]}\sum_{k=1}^{i}\,
(1+f_{k})(1+f_{h-k})\Psi_{h-k}^{[1,i]}\\
&=& w^2
\sum_{h=2}^{2N}\rho(\e_{h-i})(1+f_i)(1+f_{h-i})\Psi_{h-i}^{[1,i]}S^1_h(i)\\
&-& w^2\sum_{h=2}^{2N}\rho(\e_{h-i}) f_i f_{h-i}
\Psi_{h-i}^{[1,i]}S^2_h(i),
\end{eqnarray*}
where we have set
\begin{equation}
S^1_h(i)=\sum_{k=1}^{i}\,f_{k} f_{h-k}\Psi_{h-k}^{[1,i]},\quad
S^2_h(i)=\sum_{k=1}^{i}\,(1+f_{k})(1+ f_{h-k})\Psi_{h-k}^{[1,i]}.
\end{equation}
Now the two sums $S^1_h(i)$ and $S^2_h(i)$ are discrete
convolutions and can be evaluated for all $h$ and $i$ using the
FFT algorithm in $O(N^2\log_2N)$ operations. This can be easily
done rewriting them in the form
\begin{equation}
S_h(i)=\sum_{k=1}^{N}\,g_{k}
g_{h-k}\Psi_{h-k}^{[1,i]}\Psi_{k}^{[1,i]}, \label{eq:convd}
\end{equation}
for a suitable choice of the discrete function $g_i$.
It is well
known that for $N=2^\alpha$ with $\alpha$ integer the sum
(\ref{eq:convd}) can be computed for each $i$ via FFT in
$O(N\log_2N))=O(2^\alpha\alpha)$ operations. The total cost to
compute $S_h(i)$ for all $i$ is then $O(N^2\log_2N)$.

A better algorithm can be obtained if we rewrite the sums
$S^1_h(i)$ and $S^2_h(i)$ in the form

\begin{equation}
S_h(i)=\sum_{k=1}^{2^{\beta_i}}\,g_{k}
g_{h-k}\Psi_{h-k}^{[1,i]}\Psi_{k}^{[1,i]}, \label{eq:conv}
\end{equation}
where
\begin{equation}
\beta_i=\left\{
\begin{array}{cc}
  [[\log_2(i-1)]]+1 & \quad i>1\\
  0 &  \quad i=1
\end{array}
\right.
\end{equation}
and $[[\cdot]]$ denotes the integer part.

For each $i$ the convolution sum (\ref{eq:conv}) now can be
computed in $O(2^{\beta_i}\beta_i)$ operations.  The total cost
will be approximatively reduced by one half since $O(\sum_{i=1}^N
i\log_2i)\approx O(\frac12N^2\log_2(N))$.

Clearly once expressions $S^1_h(i)$ and $S^2_h(i)$ have been
computed the remaining two sums are of the type
\begin{equation}
g_i \sum_{h=2}^{2N} g_{h-i} \Psi_{h-i}^{[1,i]}S_h(i),
\end{equation}
which can be computed directly with $O(N^2)$ operations. Thus the
final cost for the computation of ${\tI}_1(\e_i)$ for all $i$ is
$O(N^2\log_2 N+N^2)=O(N^2\log_2N)$.

\begin{remark}
In the case of constant $\rho$ it is easy to show that expression
(\ref{eqn:qbed1f}) reduces to a double convolution sum which can
be evaluated using the FFT in only $O(N\log_2N)$ operations
instead of $O(N^2\log_2N)$.
\end{remark}

\begin{figure}[h]
\centerline{
\includegraphics[scale=0.4]{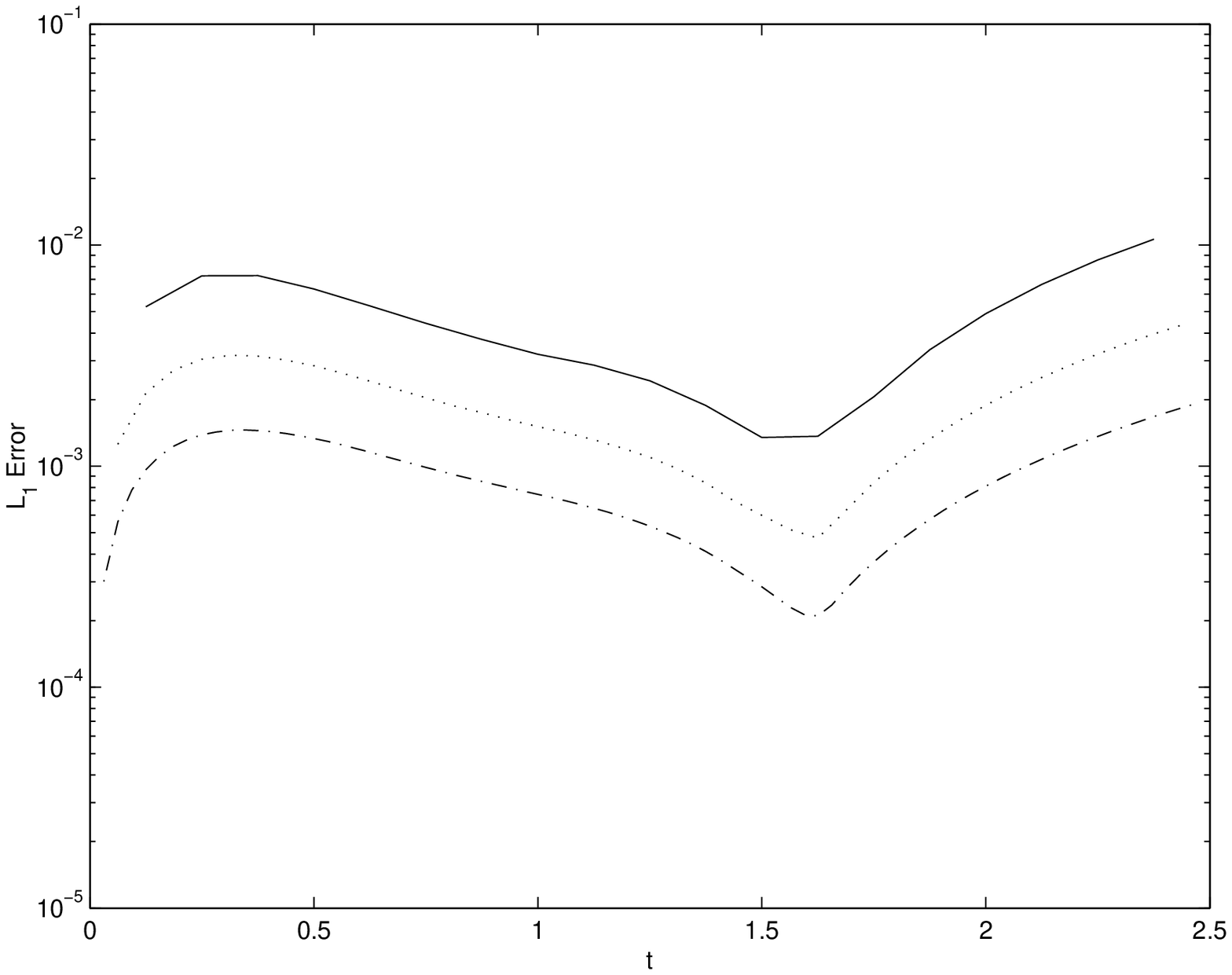}
\includegraphics[scale=0.4]{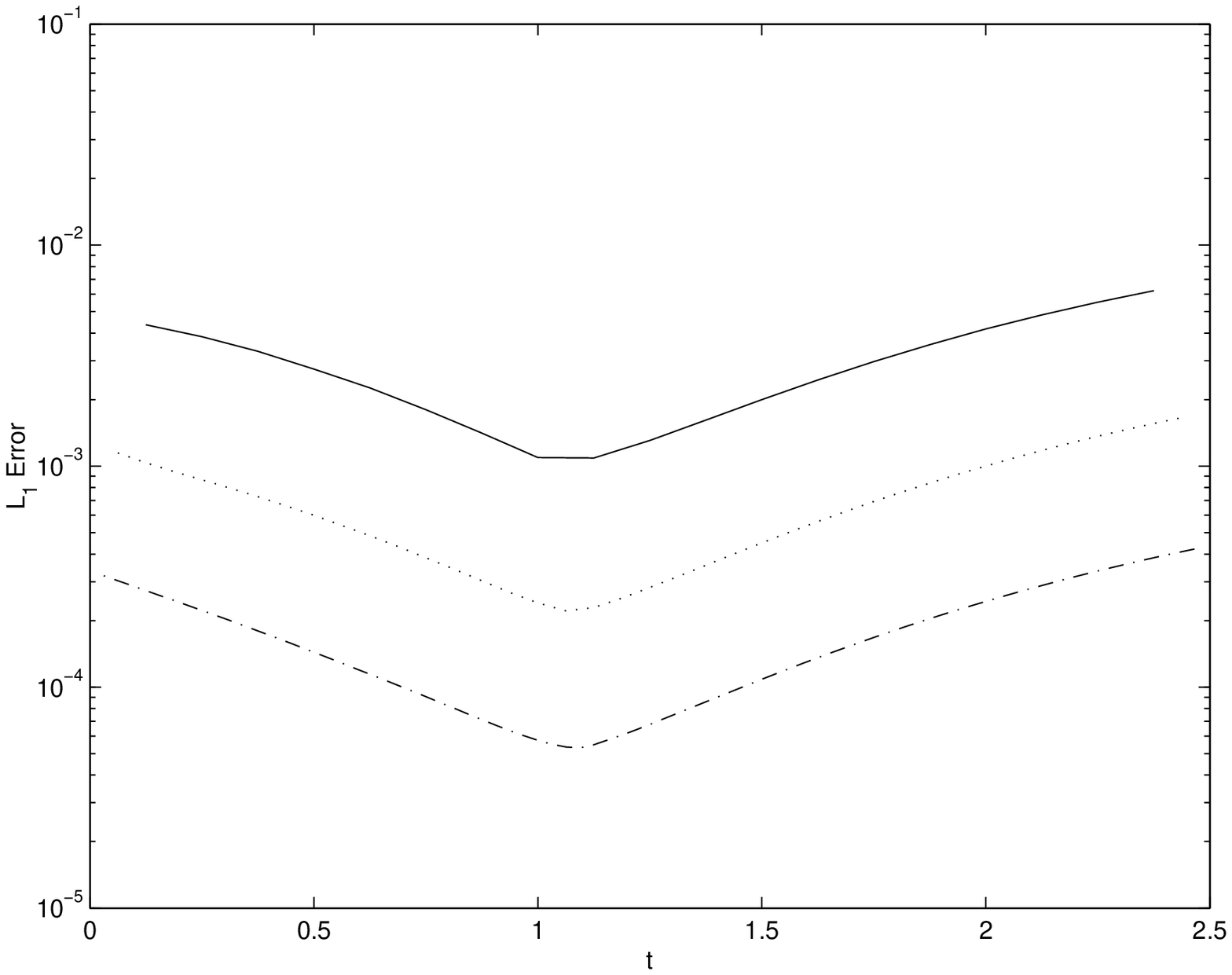}
} \caption{ \small The relative $L_1$ error for scheme QBF1 (left)
and QBF2 (right) computed with $N=20$ (solid line), $N=40$ (dotted
line), $N=80$ (dash-dot line) points for $t \in [0, 2.5]$.}
\label{fg:err1}
\end{figure}

\begin{figure}[h]
\centerline{
\includegraphics[scale=0.4]{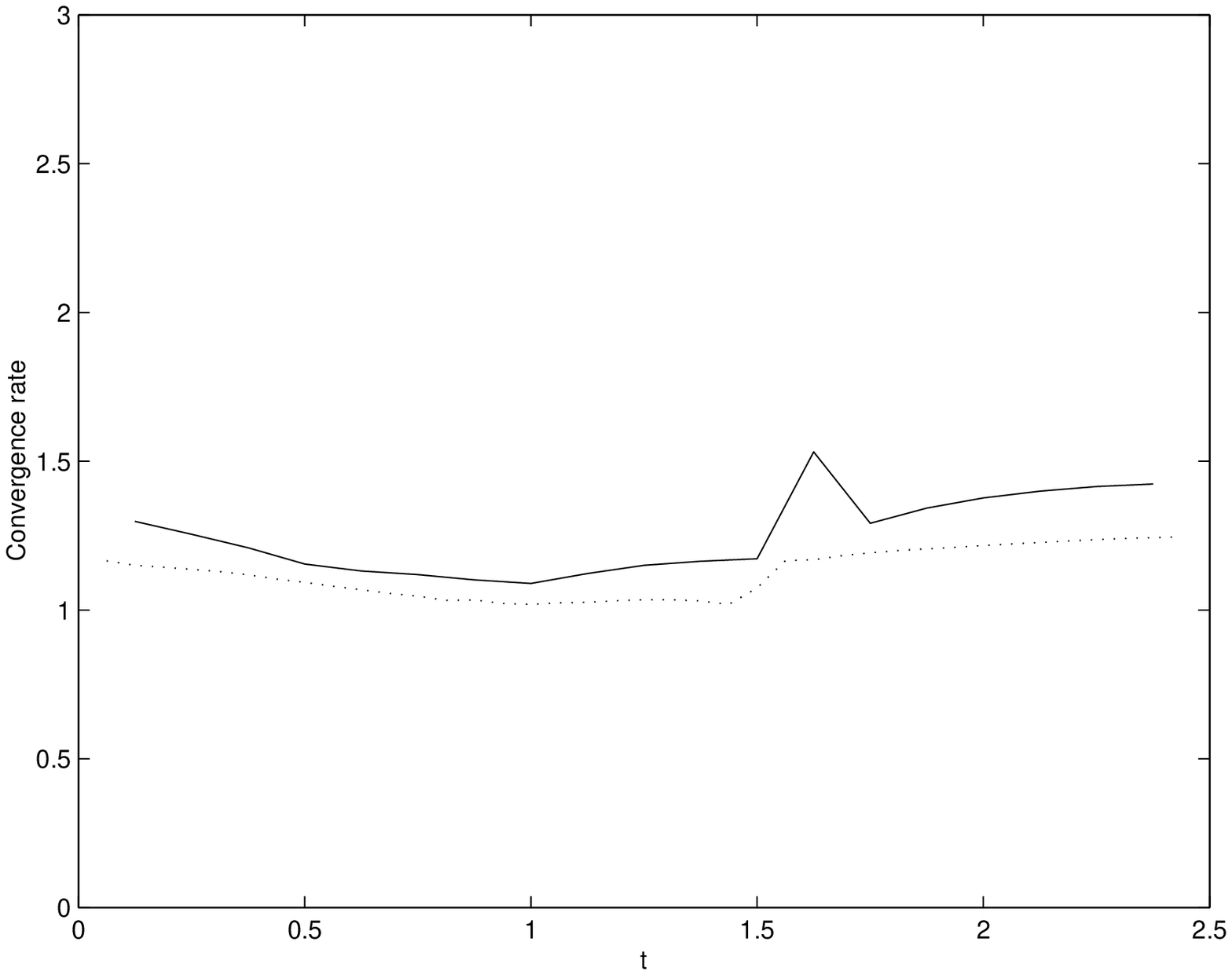}
\includegraphics[scale=0.4]{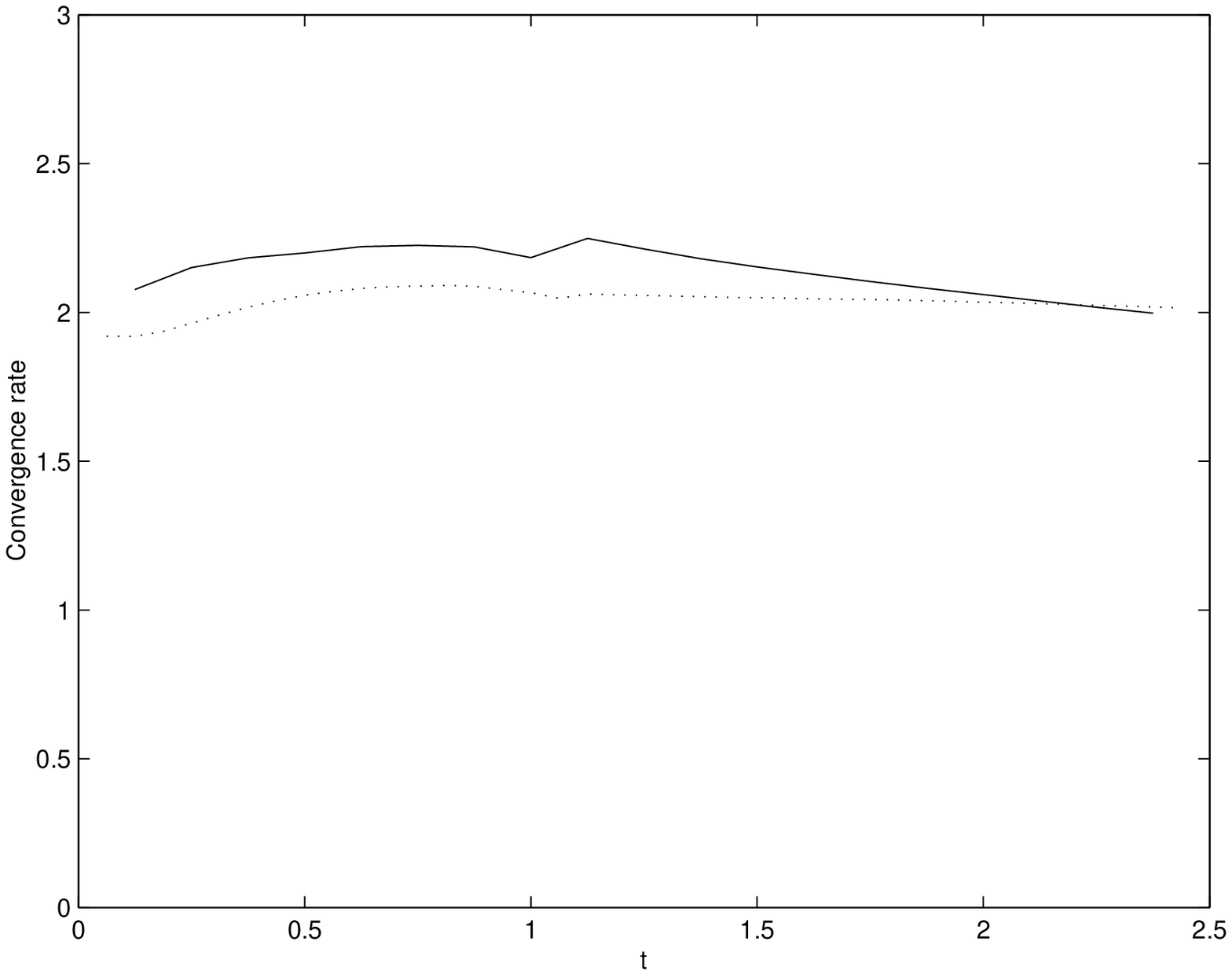}
} \caption{ \small Convergence rates for scheme QBF1 (left) and
QBF2 (right) computed with $N=20,40$ (solid line), $N=40,80$
(dotted line) points for $t \in [0, 2.5]$.} \label{fg:conv1}
\end{figure}

\section{Numerical tests and applications}
In this section we test the performance of the proposed schemes by
considering their behavior in different physical and mathematical
situations. We shall refer to the first and second order fast
schemes developed in the previous section by {\em QBF1} and {\em
QBF2} respectively. The time integration is performed with
standard first and second order explicit Runge-Kutta schemes after
dividing equation (\ref{eqn:qbord}) by $\rho(\e_i)$ and thus
rewriting the semidiscrete schemes as
\begin{eqnarray}
\nonumber \frac{\partial f_i}{\partial t}&=& w^2
\sum_{{j,l=1}\atop{1\leq k=i+j-l \leq N}}^N\,
\frac{\rho(\e_{\min})}{\rho(\e_i)} [f_k
f_l(1+f_i)(1+f_j)\\
\label{eqn:qbescheme}
\\[-.2cm]\nonumber
&-&f_i f_j (1+f_k)(1+f_l)].
\end{eqnarray}

In all our numerical tests the density of states is given by
\begin{equation}
\rho(\e)=\frac{\e^2}{2},
\end{equation}
which corresponds to an harmonic potential $V(x)$.

 Note that $0
\leq {\rho(\e_{\min})}/{\rho(\e_i)} \leq 1$ for $\e_i\neq 0$ and
that as $\e_i \to 0$ we have ${\rho(\e_{\min})}/{\rho(\e_i)} \to
1$. Furthermore since $\rho(0)=0$ the values of the distribution
function at $\e_i=0$ does not affect the discrete conservation of
mass and energy.

The schemes were implemented using the fast algorithm described in
Section \ref{fast}.

\subsection{Accuracy analysis}
The first test case has been used to check the numerical
convergence of our quadrature formulas by neglecting the time
discretization error (as usual this can be achieved either using
very small time steps or sufficiently accurate time
discretizations). The initial datum is a Gaussian profile centered
at $R/2$
\begin{equation}
f=\exp(-4(\e-R/2)^2),
\end{equation}
with $R=10$. The final integration time is $T=2.5$. We report in
Figure \ref{fg:err1} the relative errors in the $L_1-$norm
obtained with the different schemes for $N=20, 40, 80$ grid
points. As a reference solution we used the numerical result
obtained with a fine grid of $N=160$ points.

In Figure \ref{fg:conv1} the corresponding convergence rates of
the schemes are reported. As usual given two error curves $E_N$
and $E_{2N}$ corresponding to $N$ and $2N$ grid points the
convergence rate is computed as
\[
\log_2\left(\frac{E_N}{E_{2N}}\right).
\]
 The results confirm the expected
first order and second order degree of accuracy of the methods.

\begin{remark}
Since the midpoint rule, similarly to the trapezoidal rule, admits
an Euler-MacLaurin expansion we can in principle increase the
order of the method by extrapolation techniques. Unfortunately
with this approach it is difficult to keep conservations as well
as entropy inequality.
\end{remark}

\subsection{Bose-Einstein equilibrium}
\begin{figure}[h]
\centerline{
\includegraphics[scale=0.4]{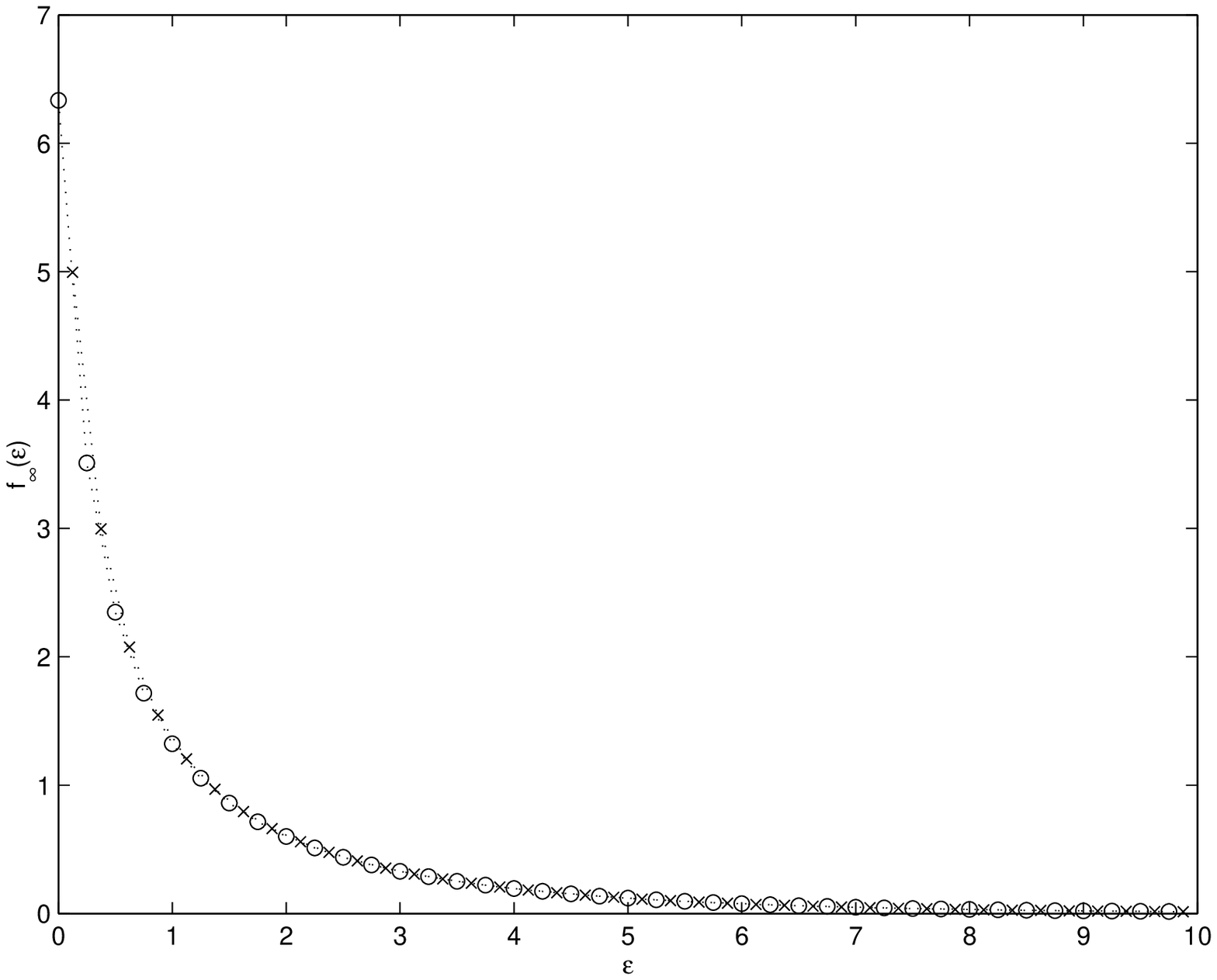}
\includegraphics[scale=0.4]{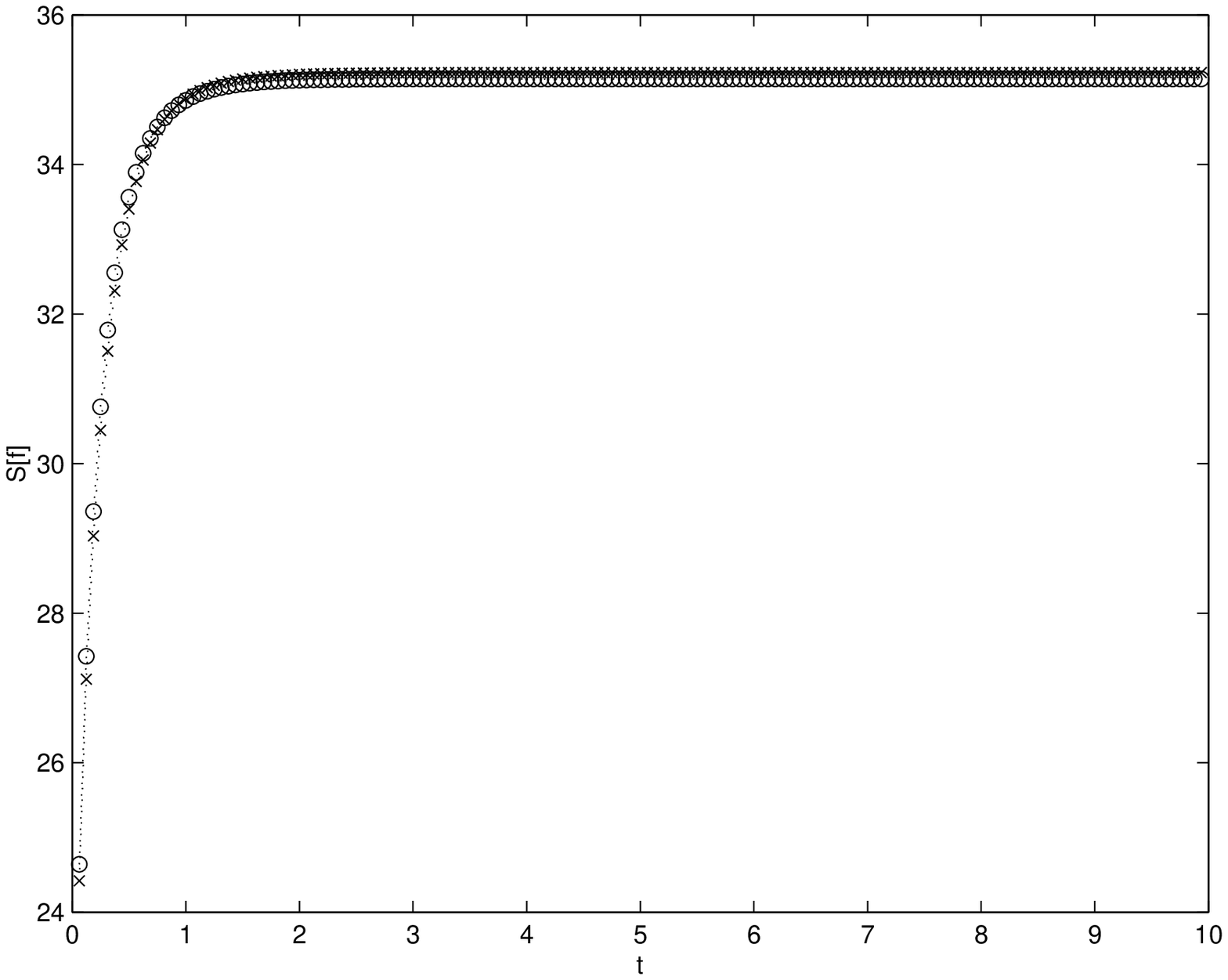}
} \caption{ \small Stationary discrete Bose-Einstein equilibrium
and entropy growth for scheme QBF1 ($\circ$) and QBF2 ($\times$)
computed with $N=40$ points.} \label{fg:Bose-Einstein}
\end{figure}

Next we consider the same initial data as in the previous section
and compute the large time behavior of the schemes for $N=40$. The
stationary solution at $t=10$ is given in Figure
\ref{fg:Bose-Einstein} for both schemes together with the
numerically computed entropy growth. As observed the methods
converge to the same stationary state given by a 'regular'
discrete Bose-Einstein distribution.

\begin{figure}[htb]
\centerline{
\includegraphics[scale=0.4]{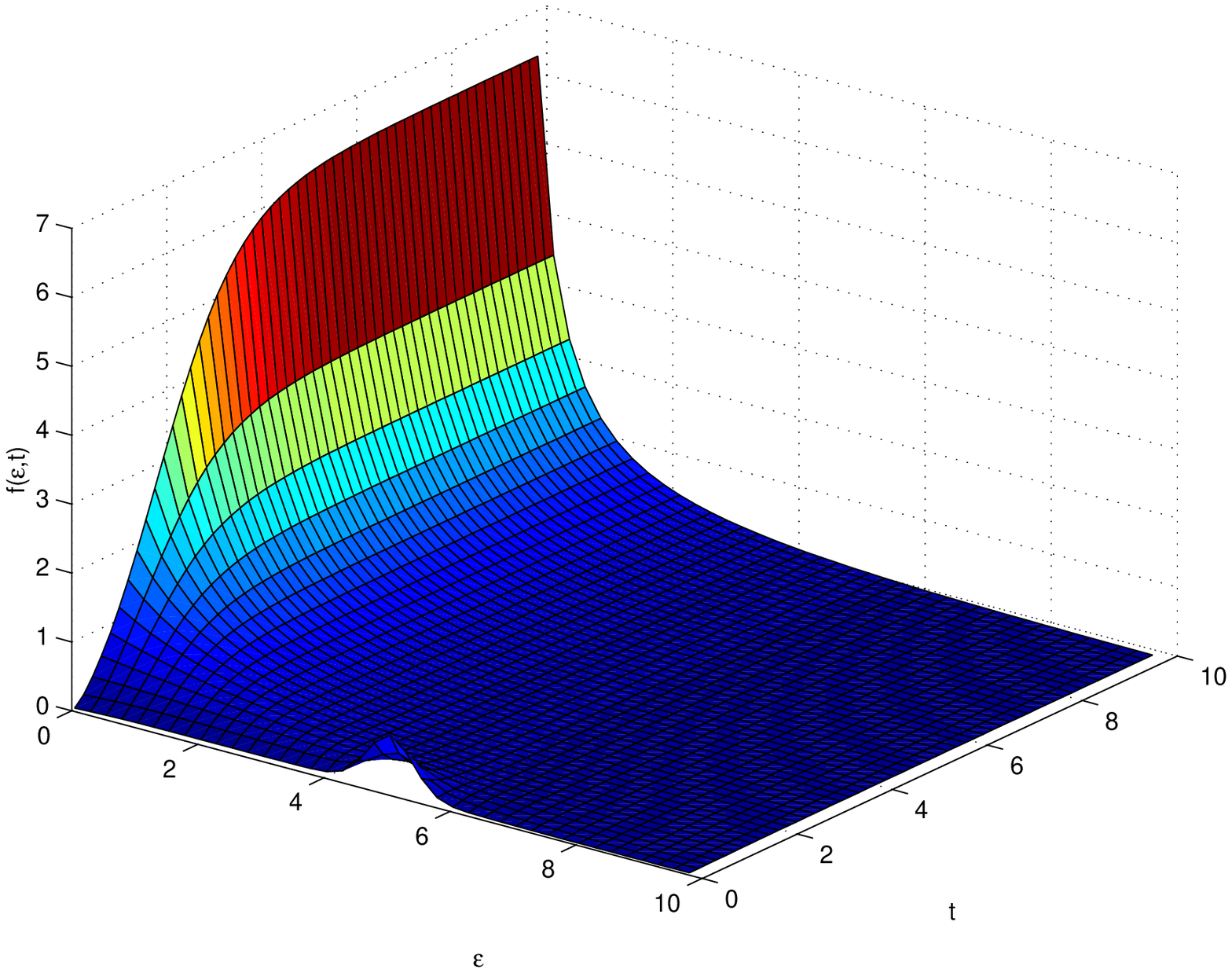}
\includegraphics[scale=0.4]{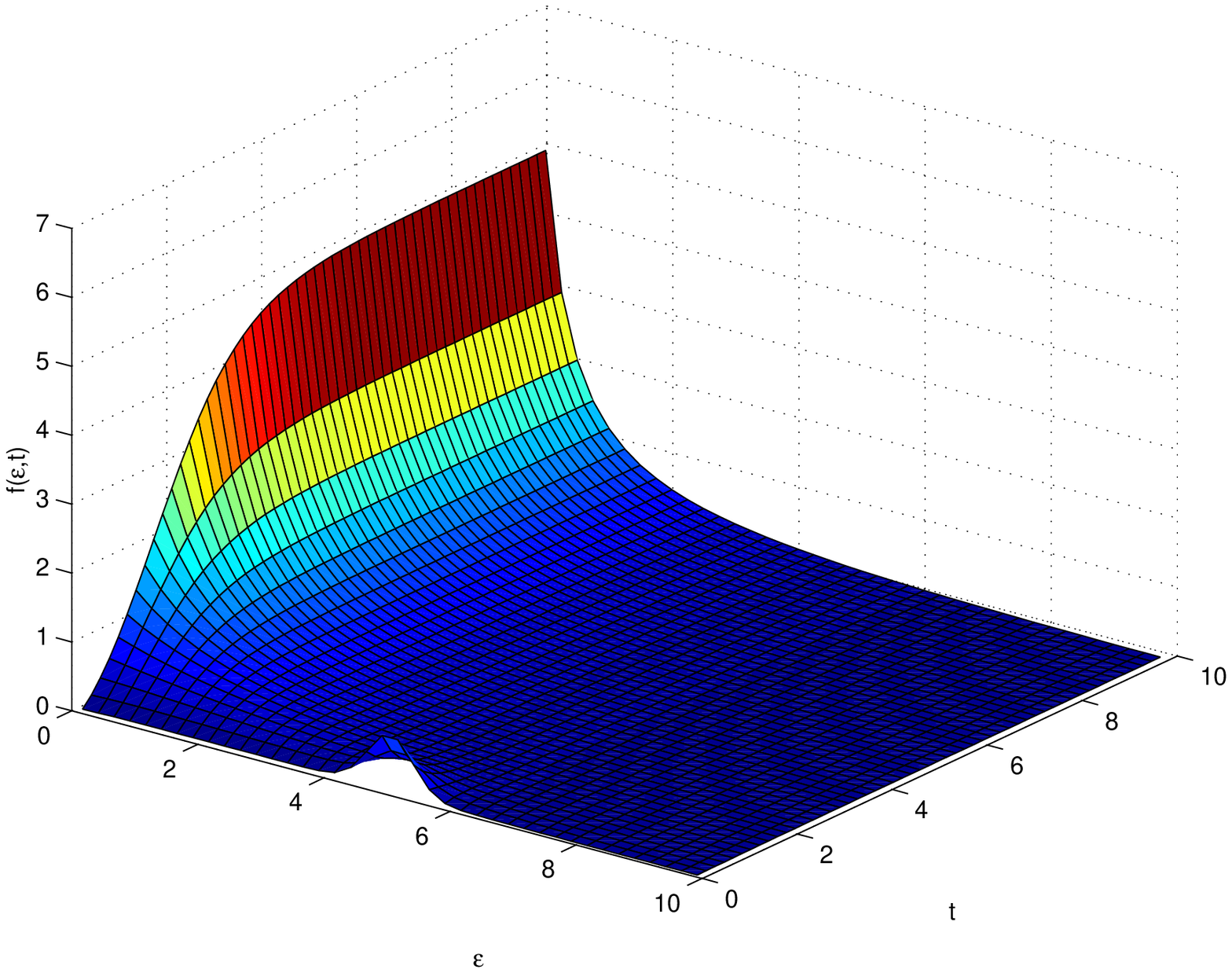}
} \caption{ \small Trend to equilibrium in time for scheme QBF1
(left) and QBF2 (right) computed with $N=40$ points.}
\label{fg:Bose-Einstein2}
\end{figure}

The trend to equilibrium in time for the two schemes is reported
in Figures \ref{fg:Bose-Einstein2}. Note that although the two
schemes agree very well there is a remarkable resolution
difference in proximity of the point $\e=0$ due to the staggered
grids of the schemes.

However since the value of $f$ at $\e=0$ does not affect the
macroscopic quantities and the entropy we can adopt a suitable
extrapolation strategy to recover a better resolution of scheme
QBF2 near $\e=0$. Since we are mostly interested in the large time
behavior of the solution we can recover the value at the zero
energy level by a steady state extrapolation. This corresponds to
assume $f$ of the form (\ref{eq:BEn}) and consequently to assign
\begin{equation}
f(0)=\frac{1}{e^{\beta}-1},\quad
\beta=\log\left(\frac{f_2+1}{f_2}\right)+\frac1{\Delta\e}\log\left(\frac{f_2(f_2+1)}{f_1(f_1+1)}\right).
\label{eq:extrap}
\end{equation}
In Table \ref{tab:one} we compare the extrapolated results at the
final computation time of scheme QBF2 for different extrapolation
methods with scheme QBF1 and with the ``exact'' steady state
solution. We remark that the values of $\alpha$ and $\beta$ for
the stationary state can be computed by inverting numerically the
equations (\ref{eq1:em})-(\ref{eq1:ee}) for $f_\infty$ given by
(\ref{eq:RBE}). The marked improvement in the resolution given by
scheme QBF2 with steady state extrapolation is evident.

\begin{table}
\begin{center}
\begin{tabular}{|c|c|c|c|c|c|}
  \hline
  Exact & QBF1&\multicolumn{4}{|c|}{QBF2 with extrapolation}\\
        \hline
        \hline
   &      & Steady state & Exponential & Cubic & Linear\\
  \hline
  7.144 & 6.335 & 7.217 & 6.449  & 6.323 & 5.994\\
  \hline
\end{tabular}
\caption{Values of $f(0)$ at $t=10$ with $N=40$ points.}
\end{center}
\label{tab:one}
\end{table}

\begin{figure}[htb]
\centerline{
\includegraphics[scale=0.4]{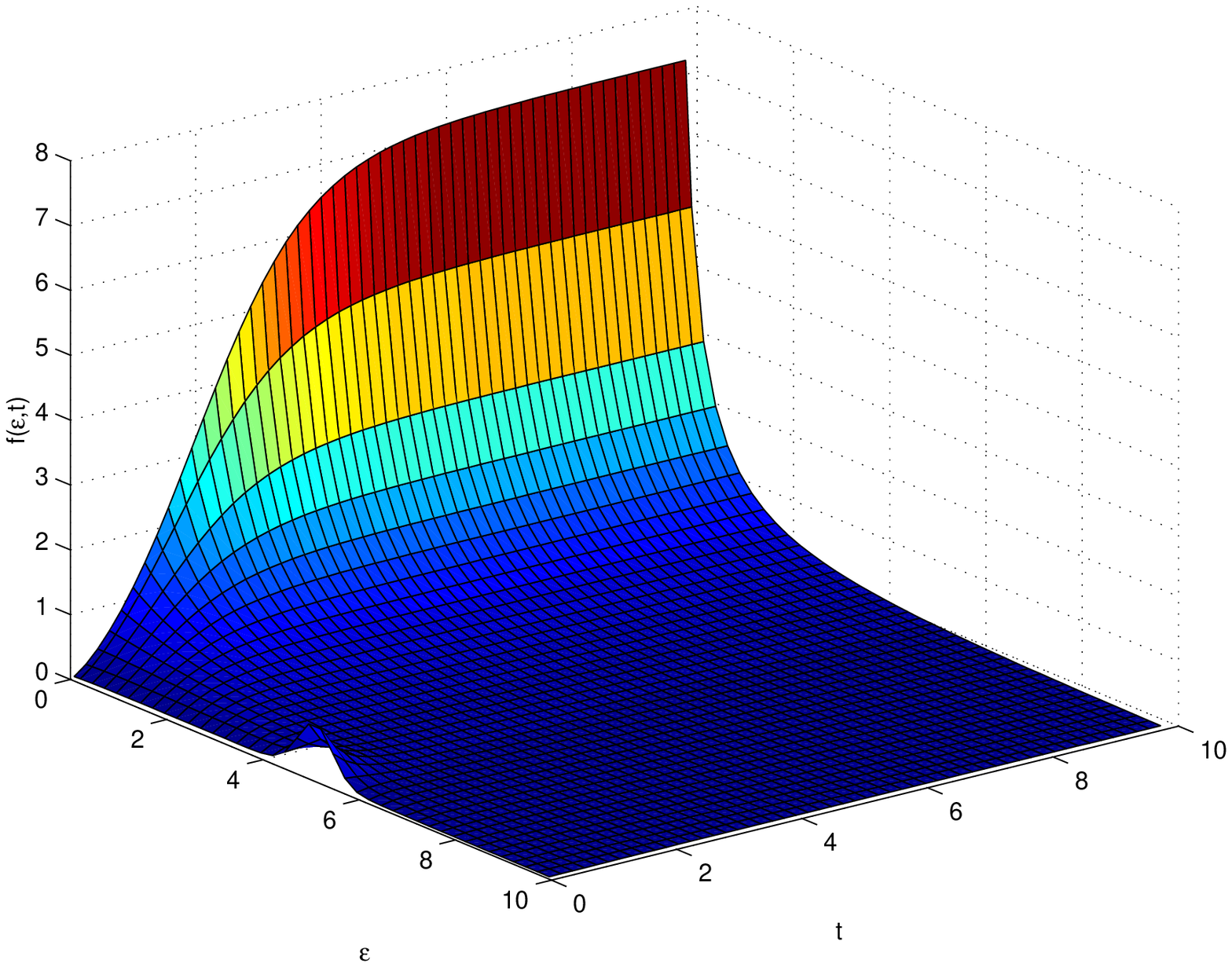}
\includegraphics[scale=0.4]{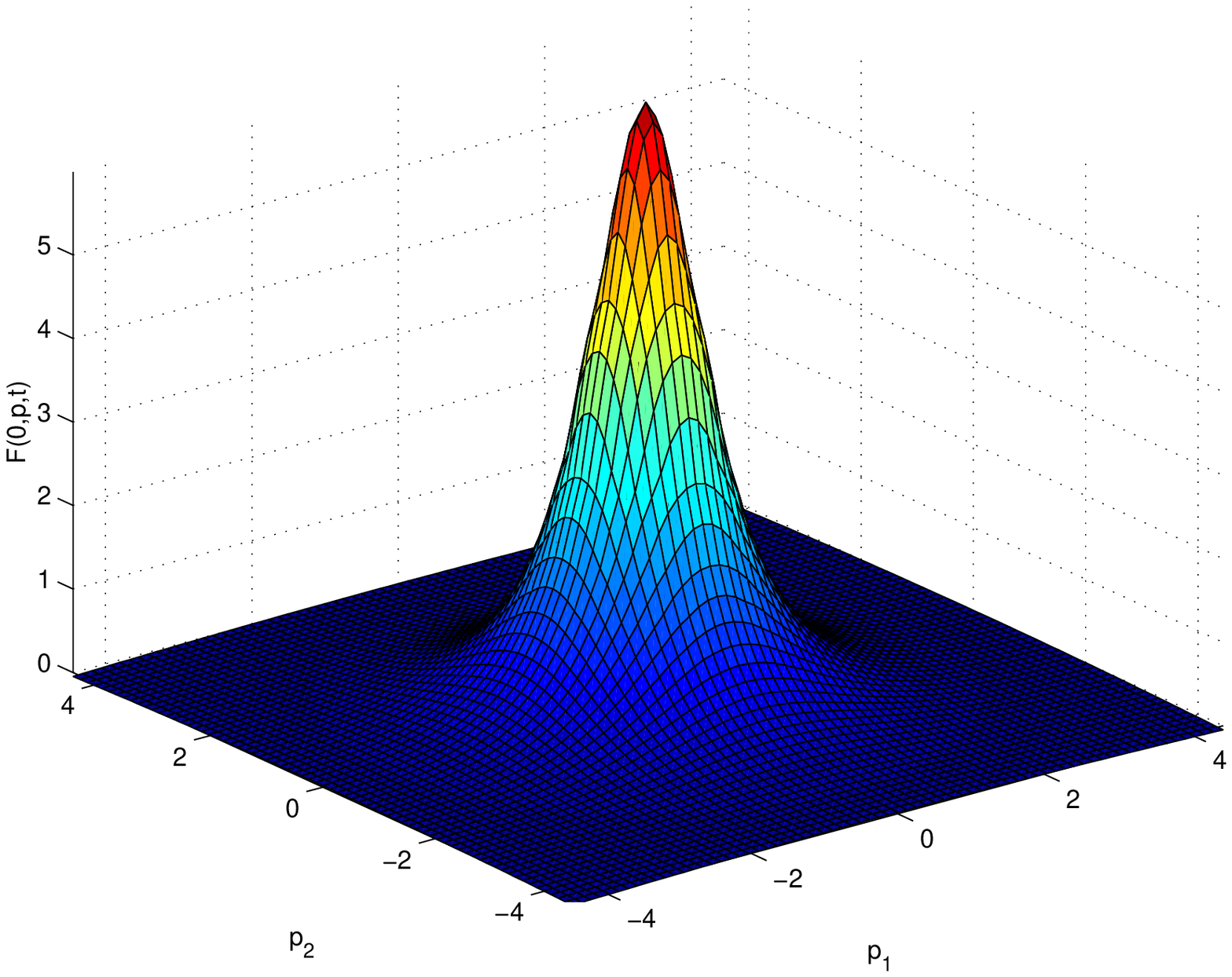}
} \caption{ \small Trend to equilibrium in time for scheme QBF2
(left) and stationary phase-space density reconstructed at $x=0$
and $p=(p_1,p_2,0)$ (right) with steady state extrapolation at
$\e=0$.} \label{fg:Bose-Einstein2b}
\end{figure}

In Figure \ref{fg:Bose-Einstein2b} we present the corresponding
result for scheme QBF2 with steady state extrapolation at $\e=0$
(as we shall always do from now on with QBF2). In the same figure
we also report the final ``steady'' solution at $t=10$ for the
phase-space density reconstructed at $x=0$ and $p=(p_1,p_2,0)$.

\subsection{Condensation}
In this test we consider the process of condensation of bosons. It
is a fundamental results of quantum statistics of bosons that
above a critical density/below a critical energy particles enter
the ground state, i.e. a Bose-Einstein condensate forms (see
\cite{GZ1},\cite{GZ2},\cite{Jaksh1},\cite{ST1},\cite{ST2}) and the
equilibrium distribution $f_\infty$ is of the form (\ref{eq:GBE})
with $\beta_-\neq 0$.


In general the evaluation of the condensate fraction as a function
of time is a challenging problem from the computational viewpoint.
If we assume the density function $f$ to be of the form
(\ref{eq:GBE}), which corresponds to the long time behavior, we
can use the following method to identify if condensation will
occur and compute the equilibrium condensate mass for a given mass
energy pair $(M,E)$.

\begin{figure}[htb]
\centerline{
\includegraphics[scale=0.5]{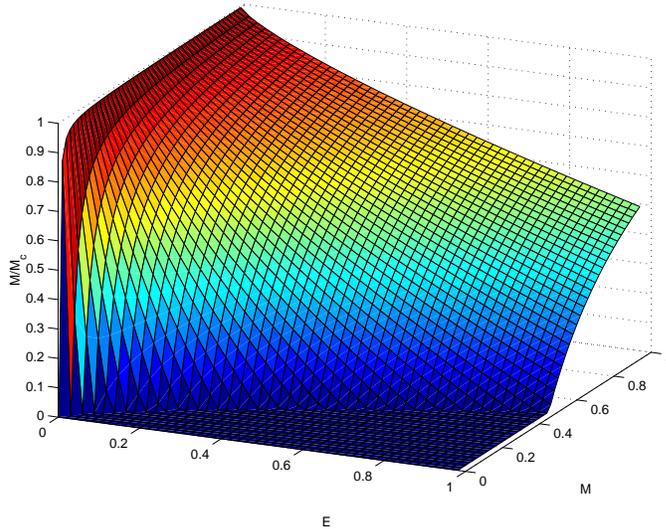}
} \caption{ \small Mass fraction of the condensate in the
mass-energy plane at the stationary state. } \label{fg:MEcond}
\end{figure}

First solve numerically for $\alpha$ the equation
\begin{equation}
E=\int_{0}^{\infty} \frac{\rho(\e)\e}{\exp(\alpha\e)-1}\,d\e.
\end{equation}
Then compute
\begin{equation}
I_\alpha=\int_{0}^{\infty} \frac{\rho(\e)}{\exp(\alpha\e)-1}\,d\e.
\end{equation}
If $I(\alpha)<M$ the mass entropy pair is critical and
condensation will take place. The condensate mass fraction in
equilibrium can then be computed
\begin{equation}
\frac{M_c}{M}=1-\frac{I_\alpha}{M}.
\end{equation}
We report in Figure \ref{fg:MEcond} the condensate mass fraction
computed with the previous method for $(M,E)\in [0,1]\times
[0,1]$.

A related challenging problem is the computation of the critical
time at which the condensate starts to form. In order to do this
we consider two different numerical indicators.

We recall that for the second order method, unlike the first order
one, due to the midpoint quadrature, we have $\e_i \neq 0$ for all
gridpoints. This makes scheme QBF2 more suitable to treat
situations where the solution is close to be singular at $\e=0$.
In particular, in such cases, it is impossible to extrapolate the
value $f(0)$ with a positive $\beta$. Thus whenever steady state
extrapolation is impossible we can assume to have formation of
condensate at $\e=0$.

For the scheme QBF1 we expect the value of $f(0)$ to increase
dramatically when formation of condensate takes place. In this
case we can use as an indicator of the formation of condensate the
expression \cite{JMP}
\begin{equation}
C_F=\frac{\Delta \e f_0}{\Delta \e \sum_{i} f_i}. \label{eq:FC}
\end{equation}


For the numerical test we choose the initial distribution in the
energy interval $[0,R]$ with $R=10$ to be\cite{ST1},\cite{ST2}
\begin{equation}
f(\epsilon)=\frac{2
\bar{f}}{\pi}\arctan(e^{\Gamma(1-\epsilon/\epsilon_0)}),
\end{equation}
with $\Gamma=5$ and $\epsilon_0=R/8$. At values of $\bar{f}$
larger than a critical $\bar{f}^*$ the formation of a condensate
occurs (see \cite{ST1},\cite{ST2} for similar results in the
homogeneous case). We choose $\bar{f}=1$, which turns on to be
supercritical. In this case the mass energy pair is
approximatively $(0.42,0.50)$ which corresponds to a condensate
mass fraction of $\approx 0.3$ at the stationary state (see Figure
\ref{fg:MEcond}).
 Using $N=320$ points and scheme QBF2 with steady
state extrapolation the condensate formation in finite time at
$t_c \approx 4.2$ is observed.

\begin{figure}[htb]
\centerline{
\includegraphics[scale=0.4]{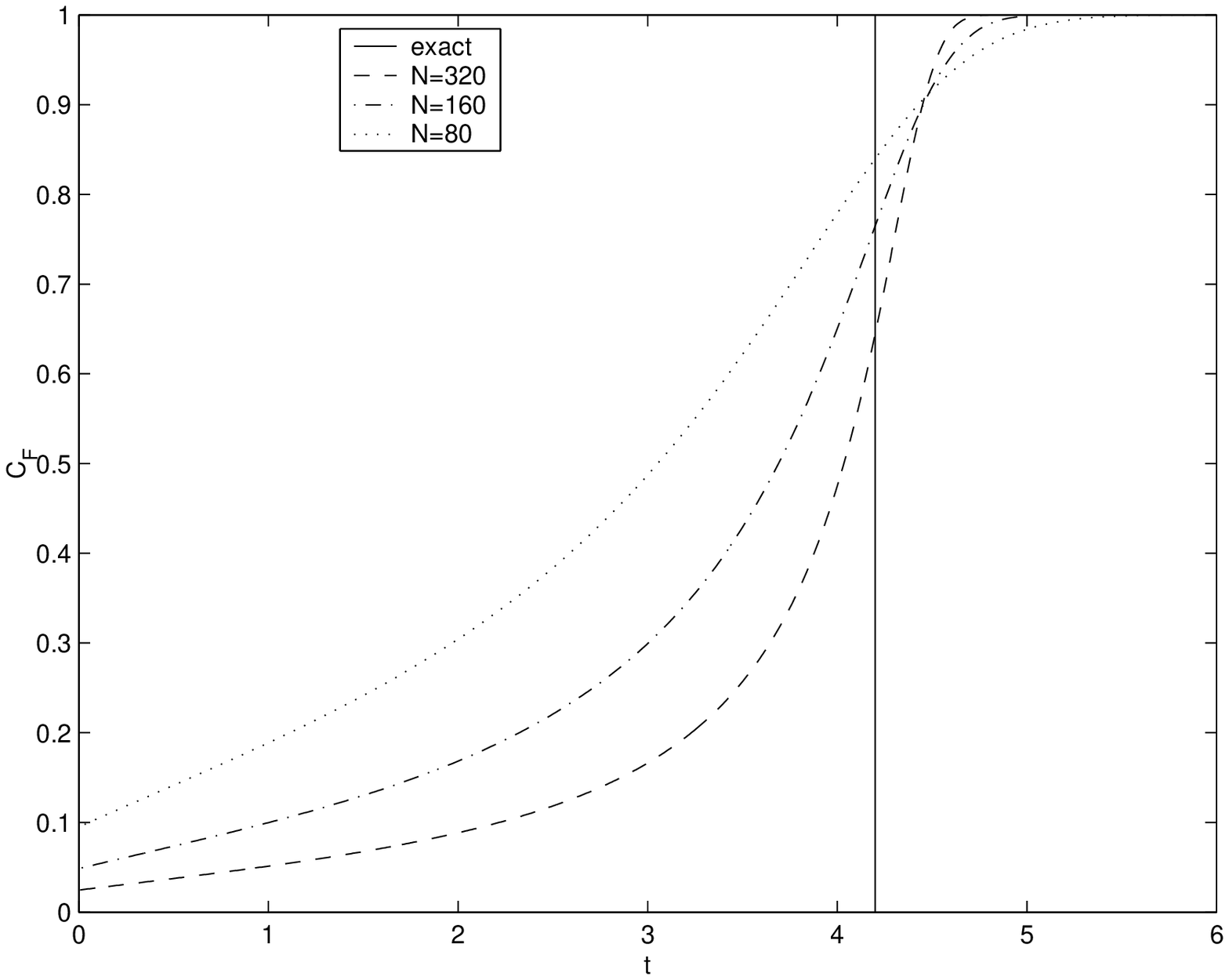}
\includegraphics[scale=0.4]{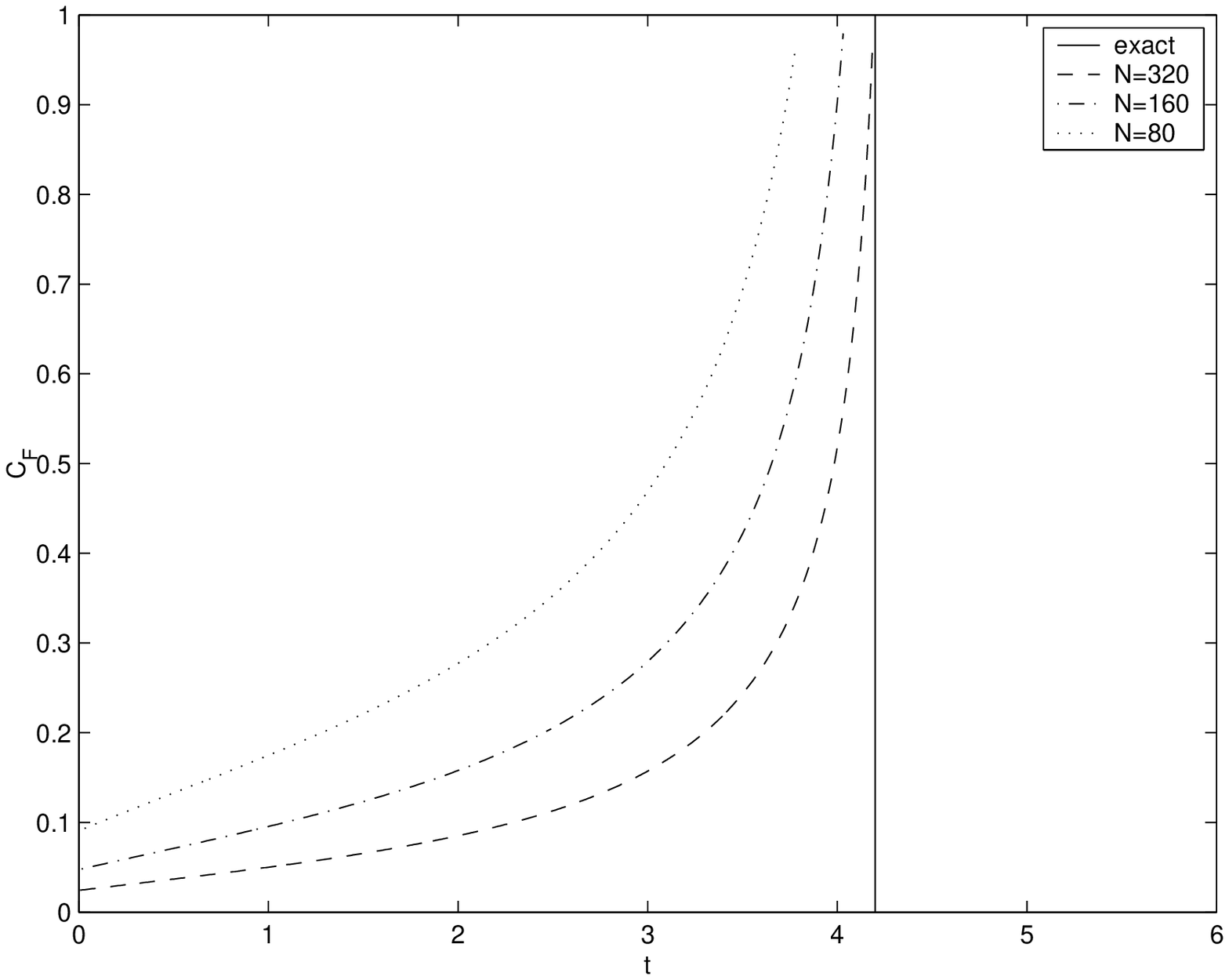}
} \caption{ \small Estimation of the critical time using the
numerical indicator $C_F$ in time for $N=80, 160, 320$ for scheme
QBF1 (left) and scheme QBF2 with steady state extrapolation
(right).} \label{fg:Bose-Einsteing}
\end{figure}

We report in Figure \ref{fg:Bose-Einsteing} the time evolution of
the indicator (\ref{eq:FC}) for scheme QBF1 and for scheme QBF2
with steady state extrapolation before the critical time. The
vertical line correspond to the critical time at which the steady
state extrapolation fails. The results indicate the numerical
convergence of the approximation (\ref{eq:FC}).

\begin{figure}[htb]
\centerline{
\includegraphics[scale=0.4]{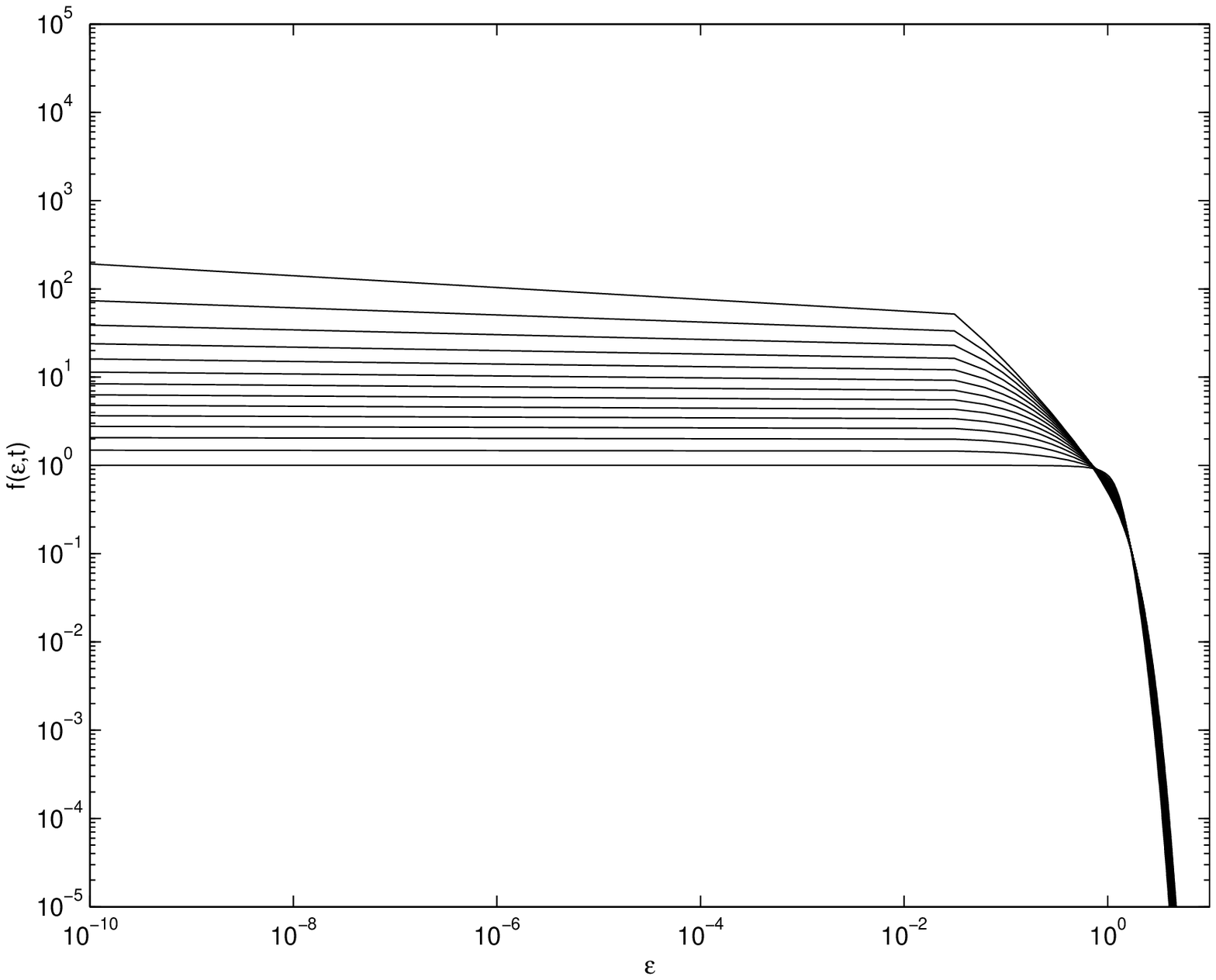}
\includegraphics[scale=0.4]{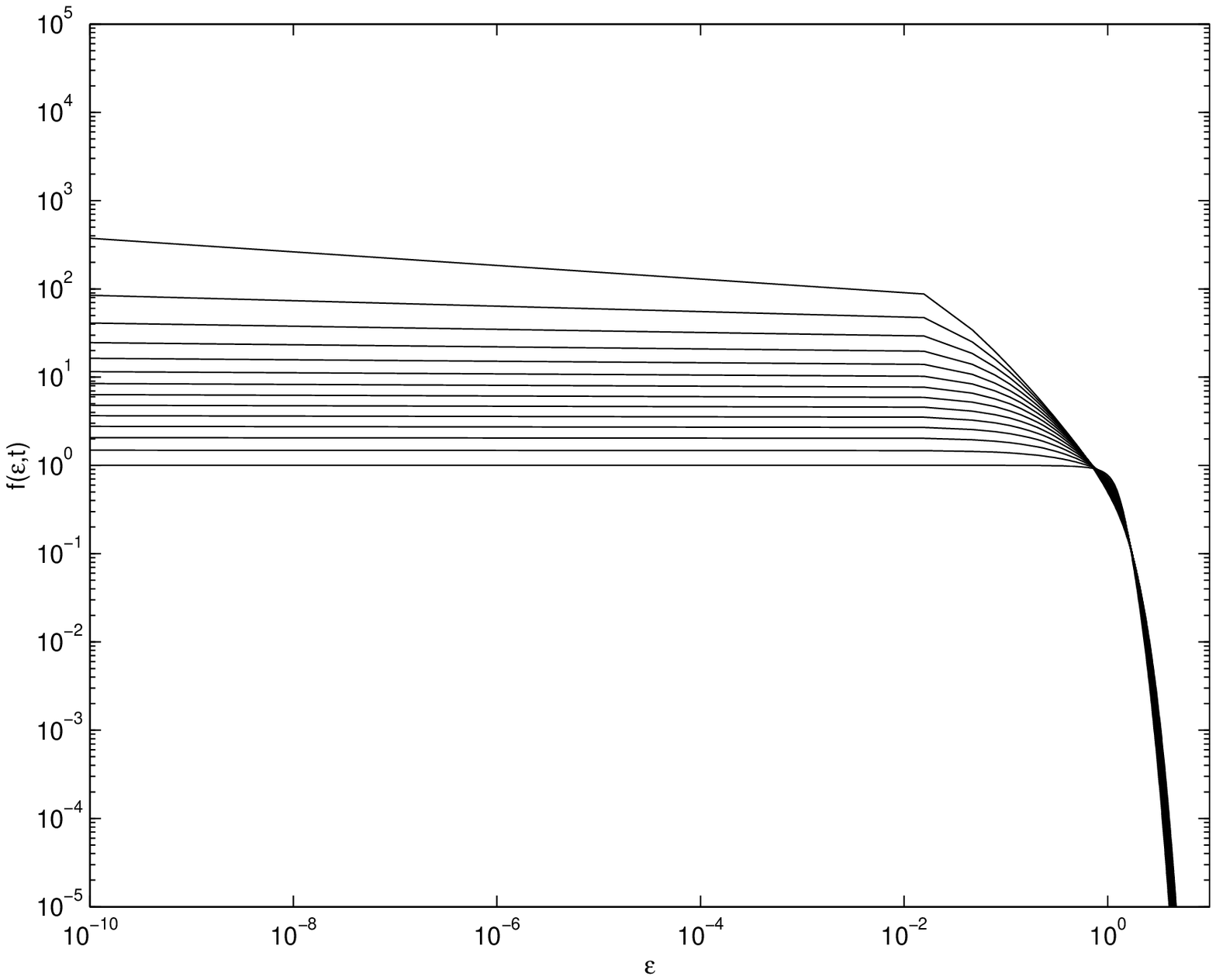}
} \caption{ \small Distribution of bosons at different times in
logarithmic scale before the critical time for scheme QBF1 (left)
and scheme QBF2 with steady state extrapolation (right) for
$N=320$.} \label{fg:conc2}
\end{figure}

\begin{figure}[htb]
\centerline{
\includegraphics[scale=0.4]{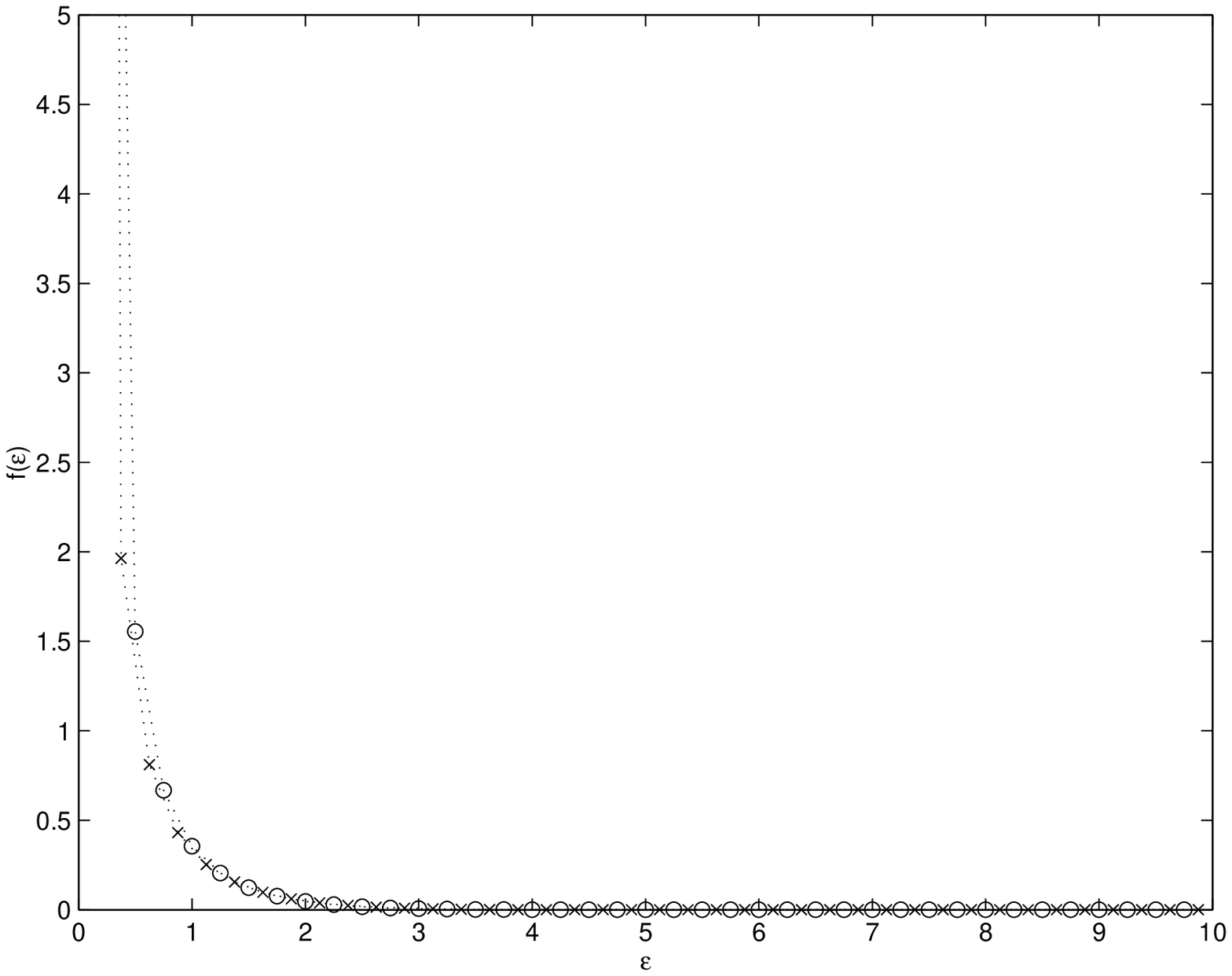}
\includegraphics[scale=0.4]{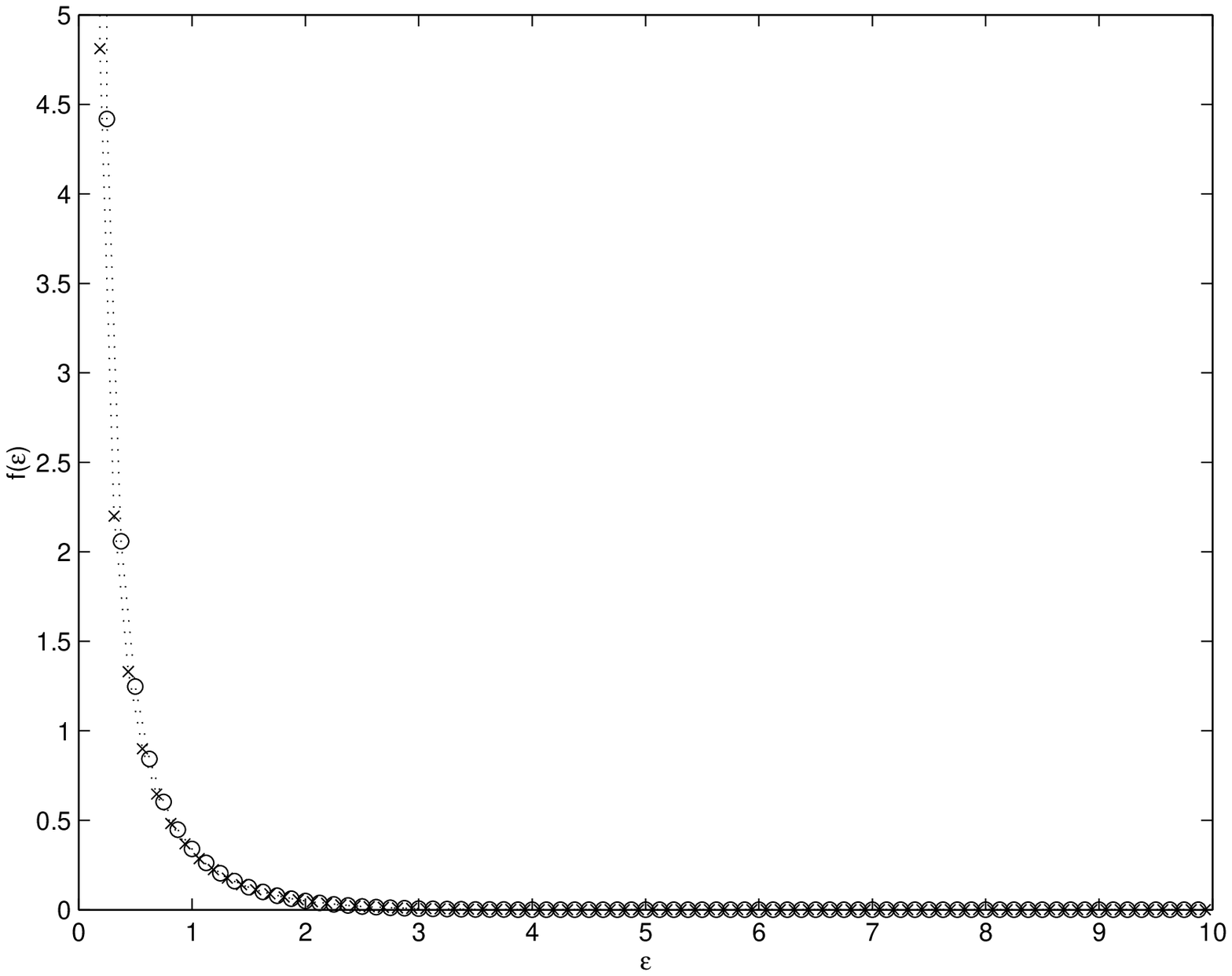}
} \caption{ \small Magnified view of the distribution of bosons
after the critical time of condensation with scheme QBF1 ($\circ$)
and scheme QBF2 ($\times$) with $N=40$ (left) and $N=80$ (right)
points at time $t=15$.} \label{fg:concm}
\end{figure}

The distribution of bosons at different times in logarithmic scale
before the critical time is shown in Figure \ref{fg:conc2} for
scheme QBF1 (left) and scheme QBF2 with steady state extrapolation
(right) .

A magnified view of the numerical solutions obtained with $N=40$
and $N=80$ points at $t=15$ shows that away from the singularity
the two schemes are still in good agreement (see Figure
\ref{fg:concm}).

\begin{figure}[htb]
\centerline{
\includegraphics[scale=0.4]{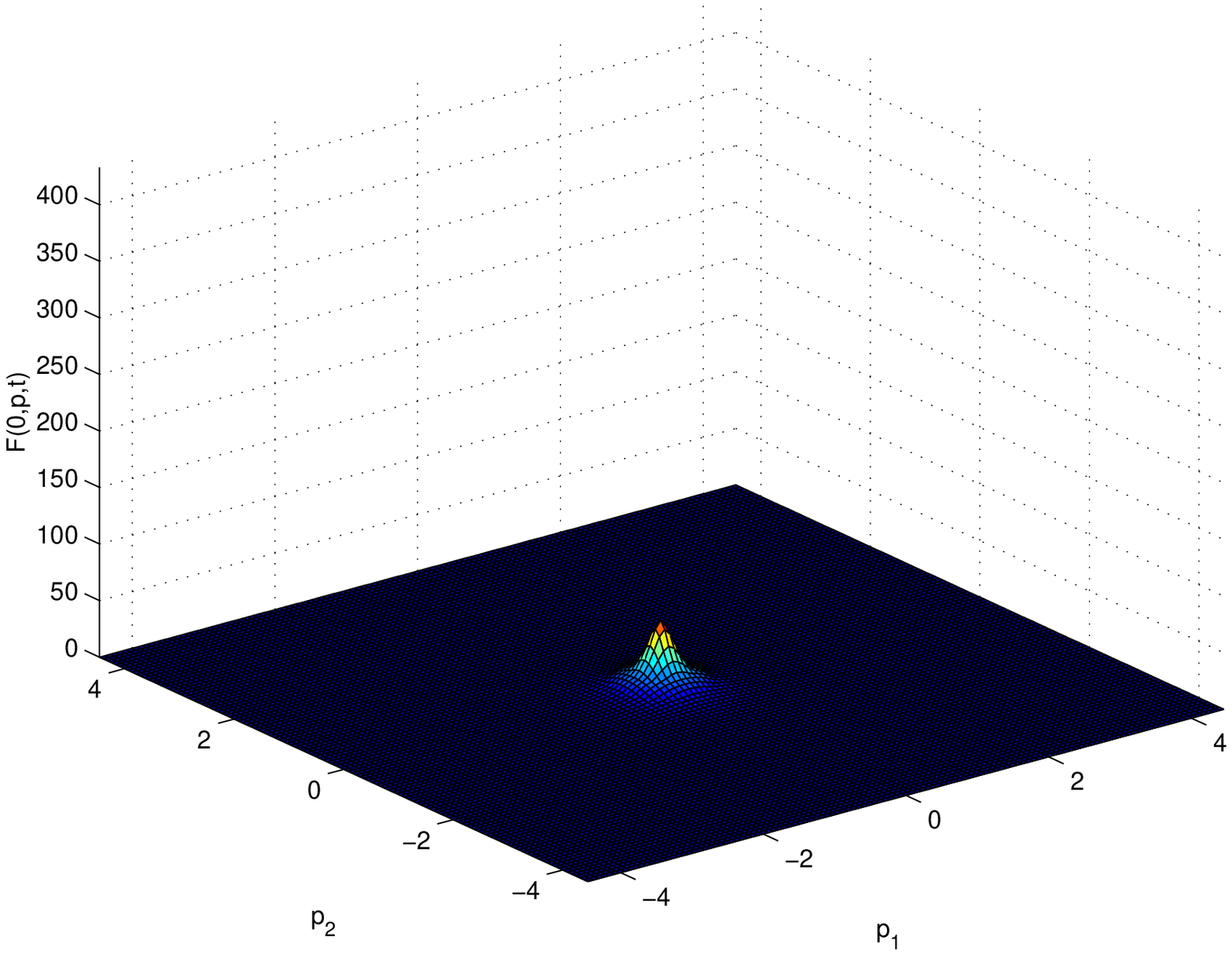}
\includegraphics[scale=0.4]{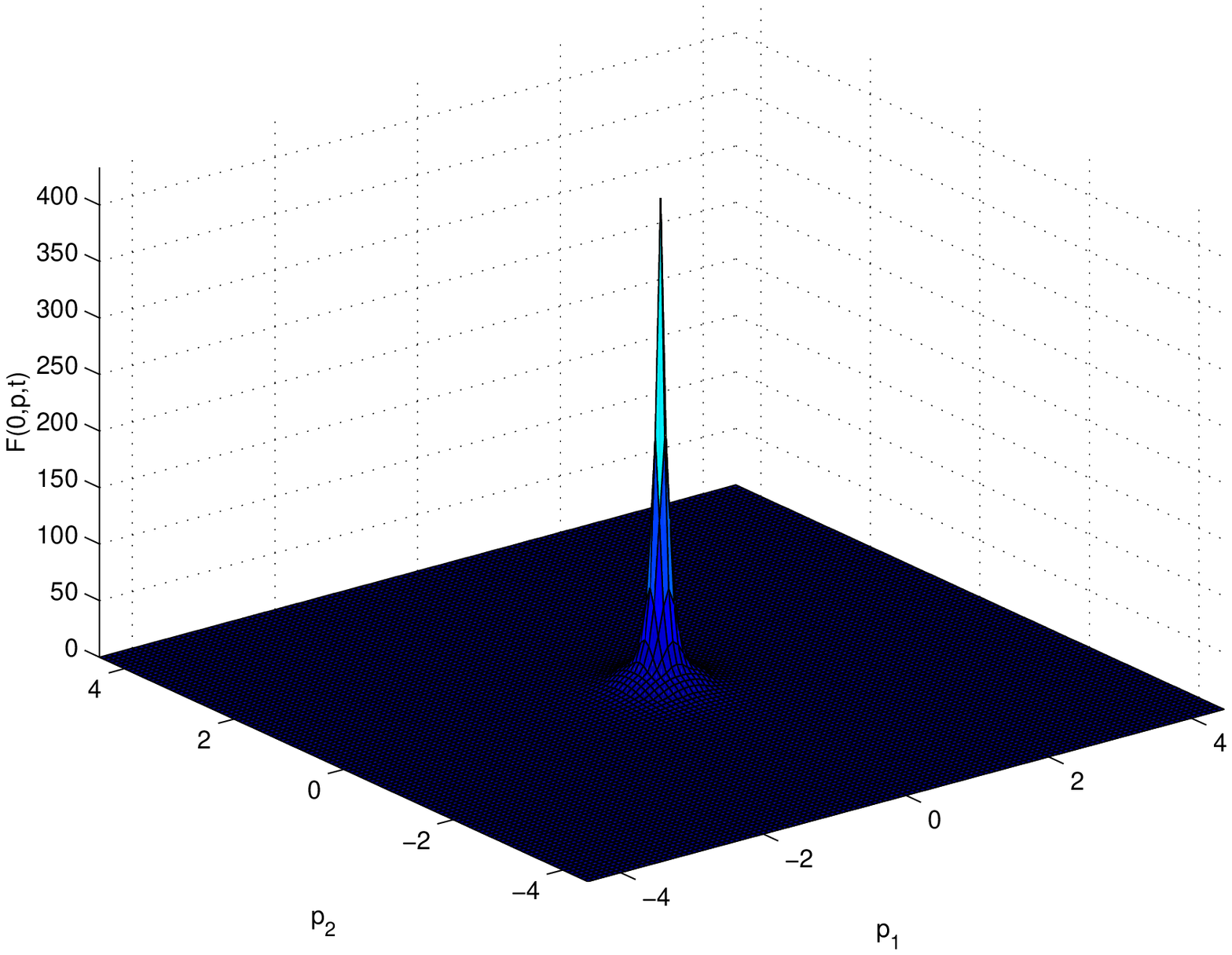}
} \caption{ \small Phase-space density reconstructed at $x=0$ and
$p=(p_1,p_2,0)$ for scheme QBF2 with steady state extrapolation at
$\e=0$ at $t=3.4$ (left) and $t=3.75$ (right).}
\label{fg:Bose-Einstein3dc}
\end{figure}

Finally in Figure \ref{fg:Bose-Einstein3dc} we also report the
phase-space density reconstructed at $x=0$ and $p=(p_1,p_2,0)$ at
two different times before the critical time. The corresponding
solution has been obtained for $N=80$ with scheme QBF2 and steady
state extrapolation.

\section{Conclusions}
We have developed first and second order fast solvers for the
Boson Boltzmann equation assuming a boson distribution which only
depends on the total energy. The methods preserve all the relevant
physical properties (conservation of mass and energy, entropy
inequality and steady states). The performance of the schemes has
been tested for both Bose-Einstein and generalized Bose-Einstein
steady states. The numerical methods have shown the capability to
describe well the challenging phenomenon of condensation of
bosons.

We remark that, to our knowledge, this is the first example of
accurate, conservative and fast deterministic numerical method for
a Boltzmann equation. Previous results were available in the
literature for Fokker-Planck-Landau type equations (see
\cite{bcdl},\cite{Lemou},\cite{PRT2}) or using some suitable
approximations of the Boltzmann equation (see the recent review
\cite{Pa} and the references therein).

Note that the present numerical methods can be applied directly
even to the case of the energy dependent quantum Boltzmann
equation for Fermions as well as the classical Boltzmann equation
of rarefied gas dynamics.

We hope to extend in the future these ideas to time dependent
potentials \cite{JMP}.

\subsection*{Acknowledgement} The authors are grateful to Dieter
Jaksch for stimulating discussions and physical explanations on
the subject of this work. We also thank the anonymous referee for
constructive suggestions.



\begin{thebibliography}{99}
\bibitem {Bo} S.N. Bose,
\newblock{Plancks Gesetz and Lichtquantenhypothese}, {\em Z. Phys.}, {\bf
26}, 178--181, (1924).

\bibitem{bcdl}
C. Buet, S. Cordier, P. Degond and M. Lemou,
\newblock{Fast algorithms for numerical, conservative, and entropy approximations of the Fokker-Planck equation},
{\it J. Comp. Phys.}, {\bf 133}, 310-322, (1997).

\bibitem{Rab} P.J.Davis, P.Rabinowitz, {\em Methods of numerical
integration}, Academic Press, (1975).

\bibitem {E1} A.Einstein,
\newblock{Quantentheorie des einatomingen idealen gases}, {\em Stiz.
Presussische Akademie der Wissenshaften Phys-math. Klasse},
Sitzungsberichte, {\bf 23}, 1--14, (1925).

\bibitem {E2} A.Einstein,
\newblock{Zur quantentheorie des idealen gases}, {\em Stiz.
Presussische Akademie der Wissenshaften Phys-math. Klasse},
Sitzungsberichte, {\bf 23}, 18--25, (1925).

\bibitem{EMV1} M.Escobedo, S.Mischler, M.A.Valle, {Homogeneous Boltzmann equation for quantum and
relativistic particles}, {\em Electron. J. Diff. Eqns.}, Monograph
04 (2003), 85 pages.

\bibitem{EMV2} M.Escobedo, S.Mischler, {On a quantum Boltzmann equation for a gas of photons},
{\em J. Math. Pures Appl.}, {\bf 9} 80, 471--515, (2001).


\bibitem{Lemou}
{M.Lemou}, \newblock{Multipole expansions for the
Fokker-Planck-Landau operator}, {\it Numerische Mathematik}, {\bf
78}, 597--618, (1998).

\bibitem{LRW}
{O.J.Luiten, M.W.Reynolds, J.T.M.Walraven}, \newblock{Kinetic
theory of evaporative cooling}, {\it Phys. Rev. A}, {\bf 53},
381--389, (1996).

\bibitem{Lu1} X.Lu, {On spatially homogeneous solutions of a modified Boltzmann equation
for Fermi-Dirac particles}, {\em J. Statist. Phys.}, {\bf 105},
353--388, (2001).

\bibitem{Lu2} X.Lu, {A modified Boltzmann equation for Bose-Einstein particles: isotropic solutions
and long-time behavior}, {\em J. Statist. Phys.}, {\bf 98},
1335--1394, (2000).

\bibitem{Pa}
L. Pareschi,
\newblock {Computational methods and fast algorithms for Boltzmann equations}, {\em Lecture Notes on the discretization
of the Boltzmann equation}, Chapter 7, Series on Advances in
Mathematics for Applied Sciences, Vol. 63, World Scientific,
(2003).

\bibitem{JMP} L.Pareschi, D.Jaksch, P.Markowich, M.Wenin,
P.Zoller, {Increasing phase-space density by varying the trap
potential in Bose-Einstein condensation}, preprint (2004).

\bibitem{PRT2}
L. Pareschi, G. Russo and G. Toscani,
\newblock Fast spectral methods for the Fokker-Planck-Landau collision operator,
\newblock {\em J. Comp. Phys}, {\bf 165}, 1--21, (2000).


\bibitem{PTV}
L. Pareschi, G.Toscani and C. Villani,
\newblock Spectral methods for the non cut-off Boltzmann equation and numerical grazing collision limit,
\newblock {\em Numerische Mathematik}, 93, pp.527-548, (2003).


\bibitem{ST1} D.V.Semikoz, I.I.Tkachev, {Kinetics of Bose
condensation}, {\em Physical Review Letters}, {\bf 74},
3093--3097, (1995).

\bibitem{ST2} D.V.Semikoz, I.I.Tkachev, {Condensation of bosons in the
kinetic regime}, {\em Physical Review D}, {\bf 55}, 489--502,
(1997).


\bibitem{Jaksh1} C.W.Gardiner, D.Jaksch, P.Zoller, {Quantum kinetic theory II}, {\em Phys.
Rev. A}, {\bf 56}, 575, (1997)

\bibitem{GZ1} C.W.Gardiner, P.Zoller, {Quantum kinetic theory}, {\em Phys. Rev. A}, {\bf 55},
2902, (1997),

\bibitem{GZ2} C.W.Gardiner, P.Zoller, {Quantum kinetic theory III}, {\em Phys. Rev. A}, {\bf 58}, 536, (1998)


\end{thebibliography}
\end{document}